\documentclass[12pt]{article}

\usepackage{mathrsfs}

\usepackage{amssymb, amsthm, amsfonts, amsxtra, amsmath}
\usepackage{latexsym}
\usepackage{mathabx}

\theoremstyle{change}

\newtheorem{Theorem}{Theorem}[section]
\newtheorem{Def}[Theorem]{Definition}
\newtheorem{Lem}[Theorem]{Lemma}
\newtheorem{Prop}[Theorem]{Proposition}

\newtheorem{Not}[Theorem]{Notation}

\date{}

\begin{document}

\hyphenation{Wo-ro-no-wicz}

\title{On a correspondence between $SU_q(2)$, $\widetilde{E}_q(2)$ and $\widetilde{SU}_q(1,1)$}
\author{Kenny De Commer\footnote{Supported in part by the ERC Advanced Grant 227458
OACFT ``Operator Algebras and Conformal Field Theory" }\\ \small Dipartimento di Matematica,  Universit\`{a} degli Studi di Roma Tor Vergata\\
\small Via della Ricerca Scientifica 1, 00133 Roma, Italy\\ \\ \small e-mail: decommer@mat.uniroma2.it}
\maketitle

\newcommand{\acnabla}{\nabla\!\!\!{^\shortmid}}
\newcommand{\undersetmin}[2]{{#1}\underset{\textrm{min}}{\otimes}{#2}}
\newcommand{\otimesud}[2]{\overset{#2}{\underset{#1}{\otimes}}}
\newcommand{\qbin}[2]{\left[ \begin{array}{c} #1 \\ #2 \end{array}\right]_{q^2}}

\newcommand{\qortc}[4]{\,\;_1\varphi_1\left(\begin{array}{c} #1  \\#2 \end{array}\mid #3,#4\right)}
\newcommand{\qortPsi}[4]{\Psi\left(\begin{array}{c} #1  \\#2 \end{array}\mid #3,#4\right)}
\newcommand{\qorta}[5]{\,\;_2\varphi_1\left(\begin{array}{cc} #1 & #2 \\ & \!\!\!\!\!\!\!\!\!\!\!#3 \end{array}\mid #4,#5\right)}
\newcommand{\qortd}[6]{\,\;_2\varphi_2\left(\begin{array}{cc} #1 & #2 \\ #3 & #4 \end{array}\mid #5,#6\right)}
\newcommand{\qortb}[7]{\,\;_3\varphi_2\left(\begin{array}{ccc} #1 & #2 & #3 \\ & \!\!\!\!\!\!\!\!#4 & \!\!\!\!\!\!\!\!#5\end{array}\mid #6,#7\right)}

\newcommand{\otimesmin}{\underset{\textrm{min}}{\otimes}}
\newcommand{\bigback}{\!\!\!\!\!\!\!\!\!\!\!\!\!\!\!\!\!\!\!\!\!\!\!\!}

\abstract{\noindent In a previous paper, we showed how one can obtain from the action of a locally compact quantum group on a type $I$-factor a possibly new locally compact quantum group. In another paper, we applied this construction method to the action of quantum $SU(2)$ on the standard Podle\'{s} sphere to obtain Woronowicz' quantum $\widetilde{E}(2)$. In this paper, we will apply this technique to the action of quantum $SU(2)$ on the quantum projective plane (whose associated von Neumann algebra is indeed a type $I$-factor). The locally compact quantum group which then comes out at the other side turns out to be the extended $SU(1,1)$ quantum group, as constructed by Koelink and Kustermans. We also show that there exists a (non-trivial) quantum groupoid which has at its corners (the duals of) the three quantum groups mentioned above.}\vspace{0.3cm}


\noindent AMS 2010 \emph{Mathematics subject classification}: 20G42; 46L65; 81R50; 16T10


\section*{Introduction}

\noindent This is part of a series of papers (\cite{DeC4},\cite{DeC5}) devoted to an intriguing correspondence between the quantizations of $SU(2)$, $\widetilde{E}(2)$ and $\widetilde{SU}(1,1)$, with the latter two groups being respectively the non-trivial two-folded covering $E(2)\rtimes \mathbb{Z}_2$ of the Euclidian transformation group of the plane, and the normalizer of $SU(1,1)$ inside $SL(2,\mathbb{C})$ (which contains $SU(1,1)$ as an index 2 normal subgroup). In a sense, their duals form a trinity of `Morita equivalent locally compact quantum groups'. There then exists a `linking quantum groupoid' combining these three quantum groups into one global structure, and it is important to understand for example the (co)representation theory of this object.\\

\noindent In this paper, we will treat the `groupoid von Neumann algebra of the linking quantum groupoid between the duals of $SU_q(2)$, $\widetilde{E}_q(2)$ and $\widetilde{SU}_q(1,1)$'. This object consists of three `corners', corresponding to linking quantum groupoids of the \emph{pairs} inside. The linking quantum groupoid between the pair consisting of the duals of $SU_q(2)$ and $\widetilde{E}_q(2)$ was treated in \cite{DeC5}. However, we will give here an alternative description which is more in line with how a second linking quantum groupoid will be presented, namely the one between the duals of $SU_q(2)$ and $\widetilde{SU}_q(1,1)$. It is this linking quantum groupoid which will be the main object of study in the present article. The third linking quantum groupoid between the duals of $\widetilde{E}_q(2)$ and $\widetilde{SU}_q(1,1)$ can then easily be obtained by a composition procedure, while the global `3$\times$3 linking quantum groupoid' is simply the three separate linking quantum groupoids pasted together.\\

\noindent Let us now describe these objects and constructions in more detail, beginning with providing some more information on the quantum groups we mentioned. We note that all $q$'s which appear in this article are real numbers satisfying $0<q<1$, and that we denote by $\mathbb{N}_0$ the set of natural numbers with $0$ \emph{ex}cluded.

\subsection*{The quantum groups $SU_q(2)$, $\widetilde{E}_q(2)$ and $\widetilde{SU}_q(1,1)$}

\noindent Of the above quantum groups, $SU_q(2)$ is the most well-known one, and the easiest to handle. It is an example of a \emph{compact quantum group} in the sense of Woronowicz (\cite{Wor2},\cite{Wor4}), and was introduced by him in \cite{Wor1} as a `twisted' or $q$-version of the ordinary $SU(2)$-group. It appears in different guises, depending on what type of functions one considers on this quantum group: polynomial, continuous or measurable. All of these viewpoints can be shown to correspond to the same `virtual object' that is $SU_q(2)$, and the passage-way between them is easy to describe. In this paper, we will only need the von Neumann algebraic picture, so we will state the definition in this context, even though this is certainly not the most suitable way to present it. We first introduce some terminology.

\begin{Def} A \emph{von Neumann bialgebra} $(M,\Delta_M)$ consists of a von Neumann algebra $M$ and a faithful normal unital $^*$-homomorphism $\Delta_M: M\rightarrow M\bar{\otimes} M$ satisfying the coassociativity condition \[(\Delta_M\otimes \iota)\Delta_M = (\iota\otimes \Delta_M)\Delta_M.\] \end{Def}

\noindent The following is a definition of $SU_q(2)$ on the von Neumann algebra level.

\begin{Def}\label{DefSUq2} Denote $I_+ = \mathbb{N}$, and denote by $\mathscr{H}_+$ the Hilbert space $l^2(I_+) \otimes l^2(\mathbb{Z})$. Consider on it the operators \[a_+=\sum_{k\in \mathbb{N}_0} \sqrt{1-q^{2k}}\,e_{k-1,k}\otimes 1,\] \[b_+= (\sum_{k\in \mathbb{N}} q^k\,e_{kk})\otimes S,\] where the $e_{ij}$ denote the standard matrix units, and where $S$ denotes the forward bilateral shift.\\

\noindent Then the \emph{von Neumann bialgebra} $(\mathscr{L}^{\infty}(SU_q(2)),\Delta_{+})$ consists of the von Neumann algebra \[\mathscr{L}^{\infty}(SU_q(2)) = B(l^2(I_+))\bar{\otimes}\mathscr{L}(\mathbb{Z}) \subseteq B(\mathscr{H}_+),\] equipped with the unique unital normal $^*$-homomorphism \[\Delta_{+}: \mathscr{L}^{\infty}(SU_q(2))\rightarrow \mathscr{L}^{\infty}(SU_q(2))\bar{\otimes}\mathscr{L}^{\infty}(SU_q(2))\] which satisfies
\[\left\{\begin{array}{l} \Delta_+(a_+) = a_+\otimes a_+ - q b_+^* \otimes b_+ \\ \Delta_+(b_+) = b_+\otimes a_+ + a_+^*\otimes b_+.\end{array}\right.\]\end{Def}

\noindent This particular von Neumann bialgebra $(\mathscr{L}^{\infty}(SU_q(2)),\Delta_+)$ will in fact have some extra structure which really qualifies it as `the space of bounded measurable functions on a (locally) compact quantum group', but we will not need this extra structure in this paper.\\

\noindent Of course, one has to verify that the above definition is meaningful. There are two ways of establishing this: the first and more natural one is to introduce first $SU_q(2)$ in a different way (by considering say its associated Hopf $^*$-algebra), and then to use its extra structure (the existence of an invariant positive state) to pass to the von Neumann algebra level, and to prove the equivalence with the above definition.\\

\noindent A second way consists of finding a unitary which implements $\Delta_+$ on the generators $a_+$ and $b_+$. The coassociativity condition will then automatically be satisfied, since it is satisfied on the generators $a_+$ and $b_+$ of $\mathscr{L}^{\infty}(SU_q(2))$. This method is a lot more computational, and makes use of some non-trivial $q$-analytic facts. However, it is this approach which is most suited for the purpose of this article.\\

\noindent Before introducing this method, let us first make some remarks on notation. We will use standard notation for all things $q$ (see \cite{Gas1}). More precisely, for $n\in \mathbb{N}\cup\{\infty\}$ and $a\in \mathbb{C}$, we denote \[(a;q)_{n} = \overset{n-1}{\underset{k=0}{\prod}} (1-q^ka),\] and \[(a_1,a_2,\ldots,a_m;q)_n = (a_1;q)_{n}(a_2;q)_n\ldots (a_m;q)_n,\] while we denote by $\,\!_r\varphi_s$ the basic hypergeometric functions. We also borrow the following notation from \cite{Koe3}.

\begin{Def}\label{DefPsi} The entire function $z\rightarrow \qortPsi{a}{b}{q}{z}$, depending on the parameters $a,b\in \mathbb{C}$, is defined as \[\qortPsi{a}{b}{q}{z} = \sum_{n=0}^{\infty} \frac{(a;q)_n(bq^n;q)_{\infty}}{(q;q)_n} (-1)^n q^{\frac{1}{2}n(n-1)}z^n.\]
\end{Def}

\noindent Then if $b\in \mathbb{C}\setminus q^{-\mathbb{N}}$, we have \[\qortPsi{a}{b}{q}{z} = (b;q)_\infty\qortc{a}{b}{q}{z}.\]

\noindent We can now state the following Proposition. We refer to Appendix \ref{ApA} for some information on how it can be deduced from well-known observations in the literature.

\begin{Prop}\label{PropBasSUq2} Writing again $I_+ = \mathbb{N}$, we denote by $P_{q^2}^{+}$ the following function on $I_+\times I_+\times I_+$: \[P_{q^2}^{+}(p,v,w) = (-q)^{p-w}q^{(p-w)(v-w)} \frac{(q^{2w+2};q^2)_{\infty}^{1/2}}{(q^2;q^2)_{\infty}(q^{2p+2},q^{2v+2};q^2)_{\infty}^{1/2}} \qortPsi{q^{2v+2}}{q^{2v-2w+2}}{q^2}{q^{2p-2w+2}},\] or equivalently, \begin{equation}\label{eqnPfunc} P_{q^{2}}^+(p,v,w) = (-1)^p q^{vw+p(v+w+1)} \frac{(q^{2v+2},q^{2p+2};q^2)_{\infty}^{1/2}}{(q^2;q^2)_{\infty}^{1/2}(q^2;q^2)_{w}^{1/2}}\qortb{q^{-2w}}{q^{-2v}}{q^{-2p}}{0}{0}{q^2}{q^2}.\end{equation}\\

\noindent Then for $r,s,t\in \mathbb{Z}$ and $p\in I_+$, the vectors \[\xi_{r,s,p,t}^+ = \underset{v-w=t}{\sum_{v,w \in I_+}} P_{q^2}^{+}(p,v,w) e_v\otimes e_{r+p-w}\otimes e_w \otimes e_{s-p+v}\] form an orthonormal basis of $\mathscr{H}_+\otimes \mathscr{H}_+ = (l^2(I_+)\otimes l^2(\mathbb{Z}))\otimes (l^2(I_+)\otimes l^2(\mathbb{Z}))$.\\

\noindent Moreover, denoting by $W_+$ the unitary \[W_+: \mathscr{H}_{+}\otimes \mathscr{H}_+\rightarrow l^2(\mathbb{Z})\otimes l^2(\mathbb{Z})\otimes \mathscr{H}_+: \xi_{r,s,p,t}^+\rightarrow e_r\otimes e_s\otimes e_p\otimes e_t,\] we have \[W_+^*(1\otimes x)W_+ = \Delta_+(x),\qquad \textrm{for all }x\in \mathscr{L}^{\infty}(SU_q(2))=B(l^2(I_+))\bar{\otimes} \mathscr{L}(\mathbb{Z}).\]
\end{Prop}

\noindent Note that there is some freedom in the choice of the $\xi_{r,s,p,t}^+$ if we only want them to implement the comultiplication. However, the above form is the most natural one to choose.\\

\noindent Let us now move on to the quantum group $\widetilde{E}_q(2)$. Also this object was introduced by Woronowicz (at least on the operator algebraic level, \cite{Wor7}), and is a $q$-version of the group of matrices \[\{\left(\begin{array}{ll} a & 0\\ b & a^{-1} \end{array}\right)\mid a,b\in \mathbb{C}, |a|=1\},\] which has an alternative abstract description as the double cover of the group $E(2)$ of Euclidian transformations of the plane. Also in this case, one has a set of different structures to consider, depending on which function algebra one is interested in. However, the passage between these structures, notably between the algebra of polynomial functions and the algebra of bounded continuous/measurable functions, is now not so straightforward as in the previous case. The main obstacles are the lack of a well-behaved invariant functional on the purely algebraic level, and the necessity to work with unbounded operators in the operator algebraic setting. This prohibits to treat the possible correspondence within a general framework. For the particular case of $\widetilde{E}_q(2)$ however, things are still quite well-behaved.\\

\noindent The following definition of $\widetilde{E}_q(2)$ is essentially the one which appears in \cite{Wor7}, but lifted to the von Neumann algebraic setting.

\begin{Def}\label{DefEq2} Denote $I_0= \mathbb{Z}$, and denote by $\mathscr{H}_0$ the Hilbert space $l^2(I_0) \otimes l^2(\mathbb{Z})$. Consider on it the unitary operator\[a_0= S^*\otimes 1,\] where $S^*$ denotes the backward bilateral shift (acting on the first factor), and the unbounded normal operator $b_0$ which has the linear span of basis vectors $e_n\otimes e_k$ as its core, with \[b_0\,e_n\otimes e_k = q^{n}\, e_n\otimes e_{k+1}, \qquad k,n\in \mathbb{Z}.\]

\noindent Then the \emph{von Neumann bialgebra} $(\mathscr{L}^{\infty}(\widetilde{E}_q(2)),\Delta_{0})$ consists of the von Neumann algebra \[\mathscr{L}^{\infty}(\widetilde{E}_q(2)) = B(l^2(I_0))\bar{\otimes}\mathscr{L}(\mathbb{Z}) \subseteq B(\mathscr{H}_0),\] equipped with the unique unital normal $^*$-homomorphism \[\Delta_{0}: \mathscr{L}^{\infty}(\widetilde{E}_q(2))\rightarrow \mathscr{L}^{\infty}(\widetilde{E}_q(2))\bar{\otimes}\mathscr{L}^{\infty}(\widetilde{E}_q(2))\] which satisfies
\[\left\{\begin{array}{l} \Delta_0(a_0) = a_0\otimes a_0 \\ \Delta_0(b_0) = b_0\otimes a_0 \dot{+} a_0^*\otimes b_0,\end{array}\right.\] where $\dot{+}$ means `the closure of the sum of two unbounded operators'.
\end{Def}

\noindent Also in this case, $(\mathscr{L}^{\infty}(\widetilde{E}_q(2)),\Delta_0)$ carries extra structure which makes $\widetilde{E}_q(2)$ eligible to be called a \emph{locally compact} quantum group.\\

\noindent As for $SU_q(2)$, it is of course not obvious on first sight that the above definition makes sense. But one can again find a unitary implementing it: the following Proposition could in principle be deduced from the results of \cite{Koe1}, but we will give another argument in the main body of the text, based on Proposition \ref{PropBasSUq2} and the results of \cite{DeC5} (see Proposition \ref{PropIdLink}).\\

\begin{Prop} \label{PropBasEq2} Writing again $I_0=\mathbb{Z}$, we denote by $P_{q^2}^{0}$ the following function on $I_0\times I_0\times I_0$: \[P_{q^2}^{0}(p,v,w) = (-q)^{p-w}q^{(p-w)(v-w)} \frac{1}{(q^2;q^2)_{\infty}} \qortPsi{0}{q^{2v-2w+2}}{q^2}{q^{2p-2w+2}}.\]

\noindent Then for $r,s,t\in \mathbb{Z}$ and $p\in I_0=\mathbb{Z}$, the vectors \[\xi_{r,s,p,t}^0 = \underset{v-w=t}{\sum_{v,w\in I_0}} P_{q^2}^{0}(p,v,w) e_{v}\otimes e_{r+p-w}\otimes e_w \otimes e_{s-p+v}\] form an orthonormal basis of $\mathscr{H}_0\otimes \mathscr{H}_0=(l^2(I_0)\otimes l^2(\mathbb{Z}))\otimes (l^2(I_0)\otimes l^2(\mathbb{Z}))$.\\

\noindent Moreover, denoting by $W_0$ the unitary \[W_0: \mathscr{H}_{0}\otimes \mathscr{H}_0\rightarrow l^2(\mathbb{Z})\otimes l^2(\mathbb{Z})\otimes \mathscr{H}_0: \xi_{r,s,p,t}^0\rightarrow e_r\otimes e_s\otimes e_p\otimes e_t,\] we have \[W_0^*(1\otimes x)W_0 = \Delta_0(x),\qquad \textrm{for all }x\in \mathscr{L}^{\infty}(\widetilde{E}_q(2))=B(l^2(I_0))\bar{\otimes}\mathscr{L}(\mathbb{Z}).\]
\end{Prop}

\noindent Finally, we have the locally compact quantum group $\widetilde{SU}_q(1,1)$ to discuss. A first question that immediately comes to mind is: why $\widetilde{SU}_q(1,1)$ and not $SU_q(1,1)$? This is because of the `no-go theorem' of Woronowicz (cf. \cite{Wor5}), which says that \emph{$SU_q(1,1)$ simply can not exist as a locally compact quantum group}. This is not as bad as it sounds: due to a key observation of Korogodsky (\cite{Kor1}), it turns out that a close companion to $SU(1,1)$ allows a $q$-deformation into a locally compact quantum group, namely the normalizer $\widetilde{SU}(1,1)$ of $SU(1,1)$ inside $SL(2,\mathbb{C})$. But one had to wait till \cite{Koe3} for the first rigorous results that this object really existed in the operator algebraic framework (and more particularly, fitted in the setting of \cite{Kus1}).\\

\noindent The presentation in \cite{Koe3} in fact \emph{started} from a concrete unitary implementing the coalgebra structure, because it turned out that the accompanying Hopf algebra structure was too weak to capture all necessary information. Hence the treatment of this quantum group had a lot more $q$-analytic machinery running in the background.\\

\noindent To state the definition of $\widetilde{SU}_q(1,1)$, we first present some auxiliary notation as in the original paper \cite{Koe3}.

\begin{Not}\label{NotSUq11} We denote $I_-^{(+)} = \mathbb{Z}$, $I_-^{(-)} = \mathbb{N}_0^- = \{m\in \mathbb{Z}\mid m<0\}$, and $I_- = I_-^{(+)}\sqcup I_-^{(-)}$, the disjoint union. We write $p\in I^{(\pm)}$ as $p_{\pm}$ when we interpret it as an element in $I_-$, and we write $\mathbf{p}$ for an indeterminate element in $\{p_{+},p_-\}$. We denote \[c: I_-\rightarrow \mathbb{Z}_2: p_{\pm}\rightarrow \pm ,\] so that $\mathbf{p} = p_{c(\mathbf{p})}$.
\end{Not}

\noindent The way in which we will present the definition is slightly different from the one in \cite{Koe3}. We again refer to appendix \ref{ApA} for more information on the equivalence between the two definitions.

\begin{Def}\label{DefSUq11} Denote by $P^{-}_{q^2}$ the following function on $I_-\times I_-\times I_-$: for $\rho,\nu,\omega\in \{\pm\}$, we put \begin{eqnarray*} P_{q^2}^{-}(p_{\rho},v_{\nu},w_{\omega}) &=& \frac{(-\rho q)^{p-w}}{\sqrt{2}}\nu^{v+1}q^{\frac{1}{2}(p+v-w)(p+v-w+1)}\\ &&\qquad \cdot \frac{(-\rho q^{-2p},-\nu q^{-2v};q^2)_{\infty}^{1/2}}{(q^4;q^4)_{\infty}(-\omega q^{-2w};q^2)_{\infty}^{1/2}} \cdot\qortPsi{-\nu q^{2v+2}}{\nu\omega q^{2v-2w+2}}{q^2}{\rho\omega q^{2p-2w+2}}.\end{eqnarray*}

\noindent Denote $\mathscr{H}_{-} = l^2(I_-)\otimes l^2(\mathbb{Z})$. Then for $\mathbf{p}\in I_-$ and $r,s,t\in \mathbb{Z}$, the vectors \[\xi_{r,s,\mathbf{p},t}^- =\underset{c(\mathbf{v})c(\mathbf{w})=c(\mathbf{p})}{\underset{v-w=t}{\sum_{\mathbf{v},\mathbf{w}\in I_-}}} P_{q^2}^-(\mathbf{p},\mathbf{v},\mathbf{w}) e_{\mathbf{v}}\otimes e_{r+p-w}\otimes e_{\mathbf{w}}\otimes e_{s-p+v} \] form an orthogonal basis of $\mathscr{H}_{-}\otimes \mathscr{H}_{-}$.\\

\noindent If we then define $W_-$ as the unitary \[W_-: \mathscr{H}_{-}\otimes \mathscr{H}_-\rightarrow l^2(\mathbb{Z})\otimes l^2(\mathbb{Z})\otimes \mathscr{H}_-: \xi_{r,s,\mathbf{p},t}^-\rightarrow e_r\otimes e_s\otimes e_{\mathbf{p}}\otimes e_t,\] the application \[x \rightarrow \Delta_-(x):=W_-^*(1\otimes x)W_-\] defines a von Neumann bi-algebra structure on $\mathscr{L}^{\infty}(\widetilde{SU}_q(1,1)) = B(l^2(I_-))\bar{\otimes}\mathscr{L}(\mathbb{Z})$.
\end{Def}

\noindent We will not present here the associated (and incomplete) Hopf algebraic picture. We refer the reader to the original paper \cite{Koe3} for this. However, it should be mentioned that the incompleteness of the Hopf algebraic picture is in some sense related to the fact that in the operator algebraic picture, certain `off-diagonal corner operators' are introduced, namely the ones intertwining the $I_-^{(+)}$ and $I^{(-)}$-part. In the Hopf-algebraic setting, there is no trace of these. For this reason, we have some doubt that $\widetilde{SU}_q(1,1)$ should really be interpreted as a $q$-deformation of $\widetilde{SU}(1,1)$. Rather, it seems to us that it is a `non-commutative blow-up' of the ordinary $SU(1,1)$-group (by changing $\mathbb{C}$ into the Morita equivalent $M_2(\mathbb{C})$, in some vague sense). For this reason, it is perhaps better to stick with Woronowicz' nomenclature `extended quantum $SU(1,1)$-group'.

\subsection*{Connecting the quantum groups by means of a linking quantum groupoid}

\noindent We now come to the notion of a \emph{linking weak von Neumann bialgebra} between these structures. A general theory of such objects was treated in \cite{DeC6} (see also \cite{DeC5} for some motivation), but we will here only present the essence of the structure for the situation at hand.\\

\noindent The observation is quite simple: consider the Hilbert spaces $\mathscr{H}_+,\mathscr{H}_0$ and $\mathscr{H}_-$ introduced in the definitions of the previous subsection, and form the direct sum Hilbert space $\mathscr{H} = \mathscr{H}_- \oplus \mathscr{H}_0\oplus \mathscr{H}_+$, which we may present in the column form $\left(\begin{array}{l}\mathscr{H}_-\\\mathscr{H}_0\\\mathscr{H}_+\end{array}\right)$. Then we have a left action on this by the direct sum von Neumann algebra \[\mathscr{L}^{\infty}(\widetilde{SU}_q(1,1))\oplus \mathscr{L}^{\infty}(\widetilde{E}_q(2))\oplus \mathscr{L}^{\infty}(SU_q(2)) = \left(\begin{array}{ccc} \mathscr{L}^{\infty}(\widetilde{SU}_q(1,1))&0&0\\ 0&\mathscr{L}^{\infty}(\widetilde{E}_q(2))&0\\0&0& \mathscr{L}^{\infty}(SU_q(2))\end{array}\right).\] But the definition of the von Neumann algebras of these quantum groups immediately suggests how this pattern can be completed at the non-diagonal entries: simply define \[\mathscr{L}(\mu,\nu) = B(l^2(I_{\mu},I_{\nu}))\bar{\otimes} \mathscr{L}(\mathbb{Z}),\] where $\mu,\nu\in \{-,0,+\}$, and with the $I_{\mu}$ defined in the definitions of the previous subsection. Then the $\mathscr{L}(\mu,\mu)$ coincide with the $\mathscr{L}^{\infty}$-von Neumann algebras of our respective quantum groups, while we now also have the action of \[ Q = \left(\begin{array}{ccc} \mathscr{L}(-,-)&\mathscr{L}(-,0)&\mathscr{L}(-,+)\\ \mathscr{L}(0,-)&\mathscr{L}(0,0)&\mathscr{L}(0,+)\\\mathscr{L}(+,-)&\mathscr{L}(+,0)& \mathscr{L}(+,+)\end{array}\right)\qquad \textrm{on}\qquad \left(\begin{array}{l}\mathscr{H}_-\\\mathscr{H}_0\\\mathscr{H}_+\end{array}\right).\] (We want to stress however that the $\mathscr{L}(\mu,\nu)$ for $\mu\neq \nu$ are \emph{not} von Neumann algebras, only Hilbert W$^*$-bimodules!)\\

\noindent The next objective is then to generalize the comultiplications of the diagonal entries $\mathscr{L}(\mu,\mu)$ to the off-diagonal parts. Also this is easy to do given all the structure at hand. To introduce this comultiplication, let us first remark that by the notation $\mathscr{L}(\mu,\nu)\bar{\otimes}\mathscr{L}(\mu,\nu)$, we mean the $\sigma$-weak closure of the algebraic tensor product of these two spaces inside $B(\mathscr{H}_{\mu}\otimes \mathscr{H}_{\mu},\mathscr{H}_{\nu}\otimes \mathscr{H}_{\nu})$. We can then collect all these tensor products together in a balanced or `$\mathbb{C}^3$-fibred' tensor product of $Q$ with itself:\[ Q* Q :=  \left(\begin{array}{ccc} \mathscr{L}(-,-)\bar{\otimes} \mathscr{L}(-,-)&\mathscr{L}(-,0)\bar{\otimes} \mathscr{L}(-,0)&\mathscr{L}(-,+)\bar{\otimes} \mathscr{L}(-,+)\\ \mathscr{L}(0,-)\bar{\otimes} \mathscr{L}(0,-)&\mathscr{L}(0,0)\bar{\otimes} \mathscr{L}(0,0)&\mathscr{L}(0,+)\bar{\otimes} \mathscr{L}(0,+)\\\mathscr{L}(+,-)\bar{\otimes} \mathscr{L}(+,-)&\mathscr{L}(+,0)\bar{\otimes} \mathscr{L}(+,0)& \mathscr{L}(+,+)\bar{\otimes} \mathscr{L}(+,+)\end{array}\right),\] which is a unital von Neumann subalgebra of $B\left(\begin{array}{l}\mathscr{H}_-\otimes \mathscr{H}_-\\\mathscr{H}_0\otimes \mathscr{H}_0\\\mathscr{H}_+\otimes \mathscr{H}_+\end{array}\right)$. The von Neumann algebra $Q*Q$ can also be identified with a corner of the tensor product $Q\bar{\otimes} Q$, namely with $h(Q\bar{\otimes} Q)h$, where $h$ is the projection \[h=1_{-}\otimes 1_- + 1_0\otimes 1_0+ 1_+\otimes 1_+,\] the $1_{\mu}$ denoting the units in $\mathscr{L}(\mu,\mu)$. (It is clear how to form then the triple fibred product $Q*Q*Q$ etc.)\\

\noindent We can now define the comultiplication maps $\Delta_{\mu\nu}$ on the off-diagonal parts: they are given by \[\Delta_{\mu\nu}: \mathscr{L}(\mu,\nu)\rightarrow \mathscr{L}(\mu,\nu)\bar{\otimes}\mathscr{L}(\mu,\nu): x \rightarrow W_{\mu}^*(1\otimes x)W_{\nu},\] where the $W_{\mu}$ were defined in the definitions of the previous subsection. Of course, one must prove that $\Delta_{\mu\nu}$ has the above range, but this is not so difficult to establish in a direct manner (see e.g.~ the proof of Proposition 3.8 in \cite{Koe3}). We can further collect these maps together into a map \[\Delta_Q:Q\rightarrow Q* Q\subseteq Q\bar{\otimes} Q\] by the formula \[ \Delta_Q(\left(\begin{array}{lll} x_{--} & x_{-0} & x_{-+} \\x_{0-} & x_{00} & x_{0+} \\ x_{+-} & x_{+0} & x_{++}\end{array}\right)) = \left(\begin{array}{lll} \Delta_{--}(x_{--}) & \Delta_{-0}(x_{-0}) & \Delta_{-+}(x_{-+}) \\ \Delta_{0-}(x_{0-}) & \Delta_{00}(x_{00}) & \Delta_{0+}(x_{0+}) \\ \Delta_{+-}(x_{+-}) & \Delta_{+0}(x_{+0}) & \Delta_{++}(x_{++})\end{array}\right),\] where $x_{\mu\nu}\in \mathscr{L}(\mu,\nu)$. This map is then obviously a faithful normal $^*$-homomorphism by definition of the maps $\Delta_{\mu\nu}$. Whether it is unital depends on the precise choice of range: if one takes $Q*Q$ as the range, then the map is unital; on the other hand, if one chooses $Q\bar{\otimes} Q$ as the range, then it is not. This is simply because the unit of $Q* Q$ is the projection $\Delta_Q(1)=h$ of $Q\bar{\otimes} Q$ which we introduced above.\\

\noindent We can now state one of the main observations in this paper.

\begin{Theorem}\label{TheoCoa} The comultiplication $\Delta_Q: Q\rightarrow Q\bar{\otimes} Q $ is coassociative.\end{Theorem}

\noindent Using the terminology of \cite{Boh1} in the von Neumann algebraic setting, this will qualify $(Q,\Delta_Q)$ as a \emph{weak von Neumann bialgebra}. Because of its particular structure, we call it a \emph{(3$\times$3-)linking weak von Neumann bialgebra} (see \cite{DeC6}). It can be interpreted as (the groupoid von Neumann algebra pertaining to) a kind of quantized groupoid with three classical objects and the duals of the quantum groups $SU_q(2)$, $\widetilde{E}_q(2)$ and $\widetilde{SU}_q(2)$ as its isotropy groups. This is also the reason why we then call these duals `Morita equivalent quantum groups', as the previous description is closely related to how Morita equivalence between (classical) group\emph{oids} is defined by means of a linking groupoid. See again \cite{DeC5} for some more intuition behind these concepts.

\subsection*{Projective corepresentations}

\noindent We now comment on the way we prove this Theorem. Our method is not straightforward, and in fact, we must admit that we have not even \emph{tried} very hard to prove Theorem \ref{TheoCoa} by direct means. This is because we hope that our method, though roundabout, is much better suited for generalization.\\

\noindent The main idea to prove Theorem \ref{TheoCoa} is the following. We first make the apparently unrelated and easy observation that for a locally compact group $G$, there is a close connection between actions on (separable) type $I$-factors on the one hand (i.e., actions on von Neumann algebras of the form $B(\mathscr{H})$ for some (separable) Hilbert space $\mathscr{H}$), and (measurable) unitary 2-cocycle functions on $G$ on the other. Indeed, given such an action, one can choose for each group element a unitary implementing the associated automorphism, and this will then provide one with an $\Omega$-projective representation for some unitary 2-cocycle function $\Omega$ (which can be taken to be measurable if the unitaries are well-chosen).\\

\noindent The philosophy is now \emph{that in the quantum setting, the proper generalization of a 2-cocycle function is a (2$\times$2-)linking weak von Neumann bialgebra}. Indeed, we showed in \cite{DeC6} that, given \emph{any} coaction of a von Neumann bialgebra on a type $I$-factor (which we then called a \emph{projective corepresentation}, see section \ref{SubProj} of the present article), one can construct from this a linking weak von Neumann bialgebra (uniquely determined up to isomorphism). Observe that we started with one von Neumann bialgebra, but that a linking weak von Neumann bialgebra has \emph{two} von Neumann bialgebras inside. Indeed, the other von Neumann bialgebra \emph{is `hidden somewhere'} in the projective corepresentation! This is a generalization of the notion of \emph{twisting} a von Neumann bialgebra by means of a unitary 2-cocycle. (In fact, our main example will arise from a genuine 2-cocycle twisting, but in a non-natural way. We will therefore not emphasize it in this paper, but refer to Proposition 4.3 of \cite{DeC6} to see the connection.)\\

\noindent The main observation then is that for $SU_q(2)$, there are \emph{two} very natural such projective representations, namely by considering the action on either the standard Podle\'{s} sphere, or on a certain $\mathbb{Z}_2$-quotient of the equatorial Podle\'{s} sphere (which can be interpreted as a quantum projective plane, \cite{Haj1}). Indeed, one can show that the von Neumann algebras associated to these quantum homogeneous spaces are both type $I$-factors. Thus one can consider their associated 2$\times$2-linking weak von Neumann bialgebras, and combine them (by a composition procedure) into a 3$\times$3-linking weak von Neumann bialgebra. This will turn out to be precisely the object described in Theorem \ref{TheoCoa}, hence proving the claimed coassociativity property in an indirect way.\\

\noindent Of course, with this discussion alone, it is not clear \emph{why} one should expect the quantum groups $\widetilde{E}_q(2)$ and $\widetilde{SU}_q(1,1)$ to pop out of these constructions. In fact, we are not sure if one can figure out a priori \emph{precisely} which quantized Lie group will appear, but one \emph{can} get some information on its associated quantized Lie algebra. Indeed, in \cite{DeC4}, an \emph{infinitesimal} picture was presented of a dual version of the object $(Q,\Delta_Q)$. This is in fact how we discovered the possibility to `deform' or `twist' $SU_q(2)$ into the other two quantum groups (see \cite{DeC4} for some more information, and for the link with actions on quantum homogeneous spaces).\\

\noindent It is our hope then that this method will allow us to obtain locally compact quantum group versions of $q$-deformations of some higher-dimensional Lie groups. (We are allowed to use the terminology \emph{locally compact quantum groups}, which correspond to von Neumann bialgebras \emph{with invariant weights}, by Proposition 3.7 of \cite{DeC6}.) Indeed, the fact that as complicated a quantum group as $\widetilde{SU}_q(1,1)$ can be obtained from this procedure, gives good hope. We want to stress that the advantage of this method is the following: actions of a compact quantum group, even on a type $I$-factor, can be described in a purely algebraic way. Then our general principle gives for free a new locally compact quantum group (and a linking structure), and the difficult analytic computations are relegated to an \emph{identification} problem, \emph{not} an existence problem. However, at the moment of writing, such generalizations have not been attempted yet, so it may well be that we are looking at an isolated phenomenon in the setting of quantum Lie groups (though it should be mentioned that on the infinitesimal level, the higher-dimensional analogues are very easily obtained).\\

\emph{Contents of the paper}\\

\noindent In \emph{the first section}, we will state the main facts concerning the theory of projective corepresentations of von Neumann bialgebras (taken from \cite{DeC6}).\\

\noindent The \emph{second section} begins with some preliminaries on the action of $SU_q(2)$ on the so-called `equatorial Podle\'{s} sphere' and on the quantum projective plane (we take \cite{Haj1} as the reference here, since both these objects are treated there together, and moreover the same conventions as ours are used). We then present the spectral decomposition of the action of $SU_q(2)$ on the quantum projective plane (but relegate the \emph{proof} to appendix \ref{ApB}), and find in this way a concrete unitary which `implements' this action.\\

\noindent In \emph{the third section}, we apply to this action the `projective corepresentation $\Rightarrow$ linking weak von Neumann bialgebra' construction we explained in the introduction, and show that the resulting object coincides with a 2$\times$2-corner of the structure described in Theorem \ref{TheoCoa}. In particular, this will show that this 2$\times$2-corner has indeed a coassociative coproduct.\\

\noindent In \emph{the fourth section}, we revisit some of the material of \cite{DeC5} to show that another of the 2$\times$2-corners of the object in Theorem \ref{TheoCoa} has a coassociative coproduct. We then end this section with the proof of Theorem \ref{TheoCoa}.\\

\noindent In \emph{Appendix \ref{ApA}}, we show that the definitions of $SU_q(2)$, $\widetilde{E}_q(2)$ and $\widetilde{SU}_q(1,1)$ we gave in the introduction are equivalent to the usual ones. In \emph{Appendix \ref{ApB}}, we carry out the computation of the spectral decomposition of $SU_q(2)$ on the quantum projective plane. In \emph{Appendix \ref{ApC}}, we prove some summation formulas for basic hypergeometric functions which were used in the article.\\

\emph{Conventions and notations}\\

\noindent By $\mathbb{N}$, we denote the set of natural numbers with zero included. By $\mathbb{N}_0$, we mean $\mathbb{N}\setminus\{0\}$. (This is important to mention since another convention is followed in \cite{Koe3}!)\\

\noindent We will mostly work with Hilbert spaces of the form $l^2(I)$, where $I$ is an index set. We then denote by $e_i$ the canonical basis vectors, and by $e_{ij}$ the corresponding canonical matrix units in $B(l^2(I))$. We denote by $\omega_{ij}$ the normal functional $\langle  e_i,\,\cdot\,e_j\rangle$ on $B(l^2(I))$ (we assume linearity in the second factor). When we write $e_{ij}$ or $e_j$ with $i$ or $j\notin I$, the element is interpreted to be zero.\\

\noindent The spatial tensor product between von Neumann algebras is denoted $\bar{\otimes}$. The ordinary tensor product between Hilbert spaces is denoted as $\otimes$. The algebraic tensor product between vector spaces is denoted as $\odot$.\\

\noindent When $A\subseteq B(\mathscr{H}_1,\mathscr{H}_2)$ and $B\subseteq B(\mathscr{H}_2,\mathscr{H}_3)$ are linear spaces of maps between certain Hilbert spaces, we will denote $B\cdot A = \{\sum_{i=1}^n b_ia_i \mid n\in \mathbb{N}_0,b_i\in B,a_i\in A\}$.\\

\noindent We use the leg numbering notation for operators on tensor products of Hilbert spaces, as is customary in quantum group theory. E.g., if $Z: \mathscr{H}^{\otimes 2}\rightarrow \mathscr{H}^{\otimes 2}$ is a certain operator, we denote by $Z_{13}$ the operator $\mathscr{H}^{\otimes 3}\rightarrow \mathscr{H}^{\otimes 3}$ acting as $Z$ on the first and third factor, and as the identity on the second factor.\\

\noindent In many formulas, we will use the notation $\pm$ and $\mp$. This means that such a formula splits up into \emph{two} formulas, one in which every $\pm$ is replaced by $+$ and $\mp$ by $-$, and one in which $\pm$ is replaced by $-$ and $\mp$ by $+$.\\

\section{Projective corepresentations}\label{SubProj}

\noindent We already used the terminology `projective corepresentation of a von Neumann bialgebra' in the introduction. Let us be a spell out the definition.

\begin{Def}\label{DefProj} Let $(M,\Delta_M)$ be a von Neumann bialgebra, and $\mathscr{H}$ a Hilbert space. By a (unitary left) \emph{projective corepresentation} of $(M,\Delta_M)$ on $\mathscr{H}$, we mean a coaction \[\alpha:B(\mathscr{H})\rightarrow M\bar{\otimes}B(\mathscr{H}),\] that is, a faithful unital normal $^*$-homomorphism satisfying the coaction property \[(\iota\otimes \alpha)\alpha = (\Delta_M\otimes \iota)\alpha.\]\end{Def}

\noindent For the applications in the subsequent sections, we will always have $\mathscr{H}= l^2(\mathbb{N})$ (in a `natural' way). \emph{For the rest of this section, we then fix a von Neumann bialgebra $(M,\Delta_M)$ and a left coaction $\alpha$ of $(M,\Delta_M)$ on $B(l^2(\mathbb{N}))$. We further assume that $M$ is represented on a Hilbert space $\mathscr{K}$ in a normal, faithful, unit-preserving way, so that we may identify $M\subseteq B(\mathscr{K})$}.\\

\noindent The following notion was introduced in \cite{DeC6}.

\begin{Def}\label{DefProjCorep} Denote $\mathscr{I} = \alpha(e_{00})(\mathscr{K}\otimes l^2(\mathbb{N}))$. The unitary \[\mathcal{G}: \mathscr{K}\otimes l^2(\mathbb{N}) \rightarrow \mathscr{I} \otimes l^2(\mathbb{N}) :\xi \rightarrow \sum_{i\in I}(\alpha(e_{0i})\xi)\otimes e_i\] will be called \emph{the implementing unitary of $\alpha$}.\end{Def}

\noindent It is easy to see that the above map $\mathcal{G}$ is indeed a well-defined unitary. Its adjoint is given by \[\mathcal{G}^*: \mathscr{I}\otimes l^2(\mathbb{N})\rightarrow \mathscr{K}\otimes l^2(\mathbb{N}): \xi\otimes \delta_i\rightarrow \alpha(e_{i0})\xi.\] For any $x\in B(l^2(\mathbb{N}))$, we then have \[\mathcal{G}^*(1\otimes x)\mathcal{G} = \alpha(x),\] which follows most easily if one takes $x$ a matrix unit. We also note that the matrix coefficients of $\mathcal{G}$ may be interpreted as \emph{Clebsch-Gordan} coefficients of $\alpha$.\\

\noindent \emph{Remark:} For the purposes of this section, it will be convenient to keep $\mathscr{I}$ the concrete Hilbert space as given above. However, in the later applications, it is more suitable to `reparametrize' $\mathscr{I}$, i.e.~ to take a unitarily equivalent copy. In this paper, this will not cause any difficulties. However, we remark that taking a different parametrization \emph{inside the same Hilbert space} has some representation-theoretic consequences (see \cite{DeC6}, Proposition 3.5).

\begin{Not}\label{NotCoob} We denote $N \subseteq B(\mathscr{K},\mathscr{I})$ for the $\sigma$-weak closure of the linear span of the set \[\{(\iota\otimes \omega_{0i})(\mathcal{G})m\mid i,j\in \mathbb{N},m\in M\}\subseteq B(\mathscr{K},\mathscr{I}),\] i.e.~ the $\sigma$-weak closure of the right $M$-module generated by the elements in the first row of $\mathcal{G}$. We denote $O = N^* \subseteq B(\mathscr{I},\mathscr{K})$ for the space of adjoints of elements in $N$. Finally, we denote by $P$ the $\sigma$-weak closure of the set $O\cdot N\subseteq B(\mathscr{I})$.\end{Not}

\noindent By definition, $N$ is a right $M$-module. It is further easy to compute that \[(\iota \otimes \omega_{0i})(\mathcal{G})^*(\iota\otimes \omega_{0k})(\mathcal{G}) = (\iota\otimes \omega_{ik})(\alpha(e_{00})) \in M, \qquad \textrm{for all }i,k\in\mathbb{N},\] so that $O\cdot N \subseteq M\subseteq B(\mathscr{K})$. We then also have that $N\cdot O \subseteq B(\mathscr{I})$ is a $^*$-algebra, and hence $P$ is a von Neumann algebra.\\

\noindent The following was proven in \cite{DeC6}, Proposition 3.6. The proof is not very hard, and follows quite immediately from the two distinguishing properties of $\mathcal{G}$, namely its unitarity and the fact that it implements $\alpha$.

\begin{Prop}\label{NotLink} Write $Q$ for the space \[\left(\begin{array}{cc} P & N \\ O & M \end{array}\right) \subseteq B(\left(\begin{array}{cc} \mathscr{I} \\ \mathscr{K}\end{array}\right)).\] Then $Q$ is a unital von Neumann subalgebra of $B(\left(\begin{array}{ll} \mathscr{I}\\\mathscr{K}\end{array}\right))$. Moreover, we have $N\cdot \mathscr{K} = \mathscr{I}$ and $O\cdot \mathscr{I} = \mathscr{K}$.\end{Prop}

\noindent The final properties imply that $(Q,e)$, with $e= \left(\begin{array}{cc} 1_P & 0 \\ 0 & 0 \end{array}\right)$, will be a \emph{linking von Neumann algebra} (between $P$ and $M$), in the sense that both $e$ and $(1-e)$ are full projections (i.e.~ $O\cdot N$ and $N\cdot O$ are $\sigma$-weakly dense in respectively $M$ and $P$). We will occasionally write the components $P,N,O,M$ as $Q_{ij}$, $i,j\in \{1,2\}$. \\

\noindent We now show that the von Neumann algebra $Q$ of the previous Proposition is endowed with more structure, namely, that it carries a coassociative comultiplication. For the proof of the following Proposition, we again refer to \cite{DeC6}, Proposition 3.6.

\begin{Prop}\label{TheoTwist} Let $\mathcal{G}$ be the unitary implementing the projective corepresentation $\alpha:B(l^2(\mathbb{N}))\rightarrow M\bar{\otimes}B(l^2(\mathbb{N}))$, and let $Q= \left(\begin{array}{cc} P & N \\ O & M \end{array}\right)$ be the von Neumann algebra as constructed above.\\

\noindent Then $\mathcal{G}\in N\bar{\otimes}B(l^2(\mathbb{N}))$, and, denoting \[Q*Q = \left(\begin{array}{cc} P\bar{\otimes}P & N\bar{\otimes}N \\ O\bar{\otimes}O & M\bar{\otimes}M \end{array}\right) \subseteq B(\left(\begin{array}{cc} \mathscr{I}\otimes\mathscr{I} \\ \mathscr{K}\otimes\mathscr{K}\end{array}\right)),\] there exists a \emph{unique} unital, normal, faithful and coassociative $^*$-homomorphism $\Gamma_Q: Q\rightarrow Q* Q$ such that $\Gamma_{Q}(Q_{ij}) \subseteq Q_{ij}\bar{\otimes}Q_{ij}$, such that the restriction of $\Gamma_Q$ to $M$ coincides with $\Delta_M$, and such that, denoting by $\Gamma_N$ the restriction of $\Gamma_{Q}$ to $N$, we have \[(\Gamma_N\otimes\iota)(\mathcal{G}) = \mathcal{G}_{13}\mathcal{G}_{23}.\]
\end{Prop}

\noindent The previous Proposition thus shows how the coaction $\alpha$ of $(M,\Delta_M)$ on $B(l^2(\mathbb{N}))$ has given rise to a linking weak von Neumann bialgebra $(Q,\Gamma_Q)$ (see Definition 0.3 of \cite{DeC6}), and in particular to a new von Neumann bialgebra $(P,\Gamma_P)=(Q_{11},\Gamma_{11})$. It is this construction method which we explained in the introduction.\\

\noindent In this paper, the projective corepresentation $\alpha$ under consideration will be the restriction of another coaction $\widetilde{\alpha}$. The following discussion is devoted to the extra structure that will be present in this case.\\

\noindent Consider the von Neumann algebra $B(l^2(\mathbb{N}))\oplus B(l^2(\mathbb{N}))$. We can identify it with $C(\mathbb{Z}_2,B(l^2(\mathbb{N})))$, the space of functions from $\mathbb{Z}_2$ to $B(l^2(\mathbb{N}))$, where we will write $\mathbb{Z}_2 = \{-,+\}$ for convenience. We can then consider the following maps: the projection maps obtained by evaluation, \[\pi_{\pm}: B(l^2(\mathbb{N}))\oplus B(l^2(\mathbb{N})) \rightarrow B(l^2(\mathbb{N})): x\rightarrow x(\pm),\] and the diagonal embedding map \[d: B(l^2(\mathbb{N}))\rightarrow B(l^2(\mathbb{N}))\oplus B(l^2(\mathbb{N})): x\rightarrow d(x)=x\oplus x.\] Then $ \pi_+ \circ d = \pi_- \circ d = \iota_{B(l^2(\mathbb{N}))}$, the identity map. We also have the flip map \[\sigma: B(l^2(\mathbb{N}))\oplus B(l^2(\mathbb{N})) \rightarrow B(l^2(\mathbb{N}))\oplus B(l^2(\mathbb{N})):x\oplus y\rightarrow y\oplus x,\] which induces an action of $\mathbb{Z}_2$. Then $d(B(l^2(\mathbb{N}))$ consists precisely of the fixed elements for $\sigma$. We will further write $e_{kl}^{(+)} = 0\oplus e_{kl}$ and $e_{kl}^{(-)} = e_{kl}\oplus 0$, and similarly for the units in these fibers: $1^{(+)} = 0\oplus 1$ and $1^{(-)} = 1\oplus 0$.\\

\noindent We will now assume that the von Neumann bialgebra $(M,\Delta_M)$ has a coaction \[\widetilde{\alpha}: B(l^2(\mathbb{N}))\oplus B(l^2(\mathbb{N})) \rightarrow M\bar{\otimes} (B(l^2(\mathbb{N}))\oplus B(l^2(\mathbb{N}))),\] which is equivariant with respect to $\sigma$: \[ \widetilde{\alpha} \circ \sigma = (\iota\otimes \sigma)\circ \widetilde{\alpha}.\] Then it is clear that $\widetilde{\alpha}$ restricts to a coaction of $M$ on $d(B(l^2(\mathbb{N})))$. We then further assume that our given coaction $\alpha$ on $B(l^2(\mathbb{N}))$ and the restriction of $\widetilde{\alpha}$ to $d(B(l^2(\mathbb{N})))$ coincide by the isomorphism $d: B(l^2(\mathbb{N}))\rightarrow d(B(l^2(\mathbb{N})))$: \[\widetilde{\alpha}\circ d  = (\iota\otimes d)\alpha.\]

\begin{Not}\label{Note} We denote \[\widetilde{\alpha}_{\pm} = (\iota\otimes \pi_{\pm})\widetilde{\alpha}: C(\mathbb{Z}_2,B(l^2(\mathbb{N})))\rightarrow M\bar{\otimes} B(l^2(\mathbb{N})).\]

\end{Not}

\noindent For the following Proposition, recall that we use the notation $\mathcal{G}$ for the projective unitary corepresentation associated to $\alpha$ (Definition \ref{DefProjCorep}), the notation $\mathscr{I}$ for the space $\alpha(e_{00})(\mathscr{K} \otimes l^2(\mathbb{N}))$, and the notation $(P,\Gamma_P)$ for the von Neumann bialgebra as constructed in Proposition \ref{TheoTwist}.

\begin{Prop}\label{PropGL} Denote \[\widetilde{e} = \widetilde{\alpha}_+(1^{(+)}-1^{(-)}) \in M\bar{\otimes} B(l^2(\mathbb{N})).\] Then there exists a self-adjoint grouplike unitary $e$ in $P$ such that $\widetilde{e} = \mathcal{G}^*(e\otimes 1)\mathcal{G}$.\end{Prop}

\noindent Recall that the group-like property means that $\Gamma_P(e) = e\otimes e$.

\begin{proof} Denote \[\breve{p}_1 = \widetilde{\alpha}_+(e_{00}^{(+)}),\]\[\breve{p}_2 = \widetilde{\alpha}_+(e_{00}^{(-)}),\] then $\breve{p}_1$ and $\breve{p}_2$ are orthogonal projections summing to $\alpha(e_{00})$. In particular, $\breve{p}_i\leq \alpha(e_{00})$, and hence they correspond to projections $p_i$ in $B(\mathscr{I})$ by the formula $\breve{p}_i = \mathcal{G}^*(p_i\otimes e_{00})\mathcal{G}$. Since $\breve{p}_i\in M\bar{\otimes} B(l^2(\mathbb{N}))$, we actually have \[p_i = \sum_{k,l\in \mathbb{N}}(\iota\otimes \omega_{0k})(\mathcal{G})(\iota\otimes\omega_{kl})( \breve{p}_i)(\iota\otimes \omega_{0l}(\mathcal{G}))^*\in P.\] Then $e := p_1-p_2$ is a self-adjoint unitary in $P$. We prove that it satisfies the conditions above.\\

\noindent Write $q_1 = \widetilde{\alpha}_+(1^{(+)})$. Then we have, for $\xi\in \mathscr{L}^2(M)\otimes l^2(\mathbb{N})$: \begin{eqnarray*}\mathcal{G}^*(p_1 \otimes 1)\mathcal{G}\xi &=& \mathcal{G}^*(p_1\otimes 1) \sum_{i} (\alpha(e_{0i})\xi)\otimes e_i \\ &=& \mathcal{G}^*(p_1\otimes 1) \sum_{i} (\widetilde{\alpha}_+(d(e_{0i}))\xi)\otimes e_i\\ &=& \mathcal{G}^*\sum_{i} (\widetilde{\alpha}_+(e_{00}^{(+)}(e_{0i}^{(-)}+ e_{0i}^{(+)}))\xi)\otimes e_i \\ &=& \mathcal{G}^*\sum_{i} (\alpha(e_{0i})\widetilde{\alpha}_+(e_{ii}^{(+)})\xi)\otimes e_i \\ &=& \sum_{i}\alpha(e_{ii})(\widetilde{\alpha}_+(e_{ii}^{(+)})\xi) \\ &=& \sum_{i}\widetilde{\alpha}_+(e_{ii}^{(+)})\xi \\ &=& q_1 \xi.\end{eqnarray*} This proves that $\widetilde{e} = \mathcal{G}^*(e\otimes 1)\mathcal{G}$. \\

\noindent We now prove that $e$ is grouplike. First, we compute that \begin{eqnarray*} (\Gamma_P(e)\otimes 1) \mathcal{G}_{13}\mathcal{G}_{23} &=& (\Gamma_P(e)\otimes 1) (\Gamma_N\otimes \iota)(\mathcal{G}) \\ &=& (\Gamma_N\otimes \iota)((e\otimes 1)\mathcal{G})\\ &=&  (\Gamma_N\otimes \iota)(\mathcal{G}\widetilde{e}) \\ &=& \mathcal{G}_{13}\mathcal{G}_{23}(\Delta_M\otimes \iota)(\widetilde{e}).\end{eqnarray*} On the other hand, \begin{eqnarray*} (e\otimes e\otimes 1)\mathcal{G}_{13}\mathcal{G}_{23} &=& ((e\otimes 1)\mathcal{G})_{13}((e\otimes 1)\mathcal{G})_{23}\\ &=& (\mathcal{G}\widetilde{e})_{13}(\mathcal{G}\widetilde{e})_{23} \\ &=& \mathcal{G}_{13}\widetilde{e}_{13}\mathcal{G}_{23}\widetilde{e}_{23}\\ &=& \mathcal{G}_{13}\mathcal{G}_{23}(\iota\otimes \alpha)(\widetilde{e})\widetilde{e}_{23}.\end{eqnarray*} So from the above two computations, we see that it is sufficient to see if \[(\iota\otimes \alpha)(\widetilde{e})\widetilde{e}_{23} = (\Delta_M\otimes \iota)(\widetilde{e}),\] i.e.~ that $\widetilde{e}$ is an $\alpha$-cocycle. Bringing $\widetilde{e}_{23}$ to the other side, and writing out the expressions with use of the coaction property $(\iota\otimes \widetilde{\alpha})\widetilde{\alpha} = (\Delta_M\otimes\iota)\widetilde{\alpha}$, this becomes, writing $f= 1^{(+)}-1^{(-)}$, \[ (\iota\otimes \alpha)(\widetilde{\alpha}_+(f)) = (\iota\otimes \widetilde{\alpha}_+)(\widetilde{\alpha}(f)(1\otimes f)).\]

\noindent Now \[(\iota\otimes \alpha)(\widetilde{\alpha}_+(f)) = (\iota\otimes \widetilde{\alpha}_+)((\iota\otimes d)\widetilde{\alpha}_+(f)), \] so it is sufficient to prove that \[(\iota\otimes d)\widetilde{\alpha}_+(f) = \widetilde{\alpha}(f)(1\otimes f),\] which is equivalent with the identity $\widetilde{\alpha}_+(f) = -\widetilde{\alpha}_-(f)$. But this follows immediately from the equivariance of $\widetilde{\alpha}$ with respect to $\sigma$, and the fact that $\sigma(f) = -f$.

\end{proof}

\noindent It is convenient to split up $\mathcal{G}$ with the aid of the projections constituting the group-like element $e$ above.

\begin{Not}\label{NotRest} Let $\widetilde{\alpha}$ be a $\sigma$-equivariant coaction of $M$ on $B(l^2(\mathbb{N}))\oplus B(l^2(\mathbb{N}))$ which restricts to $\alpha$ as above. Let $e\in P$ be the group-like element of Proposition \ref{PropGL}, and write $e=p_+-p_-$ with $p_+$ and $p_-$ orthogonal projections in $P$. Then we write \[\mathcal{G}^{(\pm)} = (p_{\pm}\otimes 1)\mathcal{G}.\]\end{Not}

\noindent It is clear that $\mathcal{G}^{(\pm)}$ are then isometries with range $(p_{\pm}\mathscr{I})\otimes l^2(\mathbb{N})$. Let us record the following fact.

\begin{Lem}\label{LemPm} Using the notation from Proposition \ref{TheoTwist} (with respect to $\alpha$) and the Notation \ref{NotRest}, we have that $\mathcal{G}^{(\pm)}\in N\bar{\otimes} B(l^2(\mathbb{N}))$ and \[(\Gamma_N\otimes \iota)\mathcal{G}^{(\pm)} = \mathcal{G}^{(+)}_{13}\mathcal{G}^{(\pm)}_{23}+\mathcal{G}^{(-)}_{13}\mathcal{G}^{(\mp)}_{23}.\]\end{Lem}

\begin{proof} As $p_{\pm}\in P$, it is immediate that $\mathcal{G}^{(\pm)} = (p_{\pm}\otimes 1)\mathcal{G}\in N\bar{\otimes}B(l^2(\mathbb{N}))$. Moreover, as $e$ is a group-like element in $P$, it follows that \[\Gamma_P(p_+) = p_+\otimes p_++p_-\otimes p_-\] and \[\Gamma_P(p_-) = p_+\otimes p_- + p_-\otimes p_+.\] As $\Gamma_N(xy) = \Gamma_P(x)\Gamma_{N}(y)$ for $x\in P$ and $y\in N$, and $(\Gamma_N\otimes \iota)\mathcal{G}=\mathcal{G}_{13}\mathcal{G}_{23}$, the formula in the statement of the Lemma follows.

\end{proof}

\section{On the action of $SU_q(2)$ on the quantum projective plane}

\subsection{The equatorial Podle\'{s} sphere and the quantum projective plane}

\noindent In \cite{DeC5}, we showed how one can apply the theory of the previous section to the action of $SU_q(2)$ on the standard Podle\'{s} sphere (i.e.~ the one which arises as the quotient space by the $S^1$-action). In this paper, we will need to consider another Podle\'{s} sphere, namely the \emph{equatorial one} (which is the other extreme point in the moduli space of Podle\'{s} spheres). Let us recall the definition in the version which will be of most use to us. The equivalence of this definition with the ordinary one as a universal object can be found for example in \cite{Haj1}. We also refer to that paper for the notion of the \emph{quantum projective plane} (note that \emph{this} `projective' is unrelated to the one of the previous section!). We will keep using the notations $s,d$ and $\pi_{\pm}$ we introduced near the end of the previous section.

\begin{Def} Denote by $Y^{(\pm)}$ and $W^{(\pm)}$ the following operators on $l^2(\mathbb{N})$: \[Y^{(\pm)} = \pm \,\sum_{k\in \mathbb{N}_0} \sqrt{1-q^{4k}}\; e_{k-1,k},\]\[W^{(\pm)}\rightarrow \pm\,\sum_{k\in \mathbb{N}} q^{2k}\, e_{kk}.\] Consider then $Y,W\in C(\mathbb{Z}_2,B(l^2(\mathbb{N}))) = B(l^2(\mathbb{N}))\oplus B(l^2(\mathbb{N}))$, with \[ Y(\mu):= Y^{(\mu)},\]\[W(\mu) := W^{(\mu)}.\]

\noindent Then the unital C$^*$-algebra $C(S_{q\infty}^2)$ generated by $Y$ and $W$ is called the space of continuous functions on the \emph{equatorial Podle\'{s} sphere} $S_{q\infty}^2$.\\

\noindent The $\mathbb{Z}_2$-action $\sigma$ on $C(\mathbb{Z}_2,B(l^2(\mathbb{N})))$ restricts to an action of $\mathbb{Z}_2$ on $C(S_{q\infty}^2)$, given on the generators as $\sigma(Y)=-Y$ and $\sigma(W)=-W$. We denote the space of $\mathbb{Z}_2$-fixed elements as $C(\mathbb{R}P_q^2)$, and call it the space of continuous functions on the \emph{quantum projective plane} $\mathbb{R}P_q^2$.
\end{Def}

\noindent It is well-known that $SU_q(2)$ has a natural ergodic action $\widetilde{\alpha}$ on $S_{q\infty}^2$. Since it commutes with the $\mathbb{Z}_2$-action (see e.g.~ Remark 4.2 in \cite{Haj1}), we then also have an ergodic action $\alpha$ of $SU_q(2)$ on $\mathbb{R}P_q(2)$. Now in Lemma 6.5 of \cite{Tom1}, it is shown that the ensuing $SU_q(2)$-invariant state $\widetilde{\omega}$ on $C(S_{q\infty}^2)$ is obtained by applying the functional \[\widetilde{\omega}(x) = \frac{1}{2}(\omega(x(+)) + \omega(x(-)))\] on $C(S_{q\infty}^2)\subseteq C(\mathbb{Z}_2,B(l^2(\mathbb{N})))$, where $\omega$ is $\textrm{Tr}(\,\cdot\, D)$ with $D$ the trace class operator $(1-q^2)\textrm{Diag}(q^{2k})$. From this, and the fact that $C(S_{q\infty}^2)'' = C(\mathbb{Z}_2,B(l^2(\mathbb{N})))$, it follows that the von Neumann algebra $\mathscr{L}^{\infty}(S_{q\infty}^2) = \pi_{\widetilde{\omega}}(C(S_{q\infty}^2))''$, with $\pi_{\widetilde{\omega}}$ the GNS-representation, may be identified with $C(\mathbb{Z}_2,B(l^2(\mathbb{N})))$, in such a way that $\pi_{\widetilde{\omega}}(Y)$ and $\pi_{\widetilde{\omega}}(W)$ coincide with respectively $Y$ and $W$.\\

\noindent It then also follows that $\mathscr{L}^{\infty}(\mathbb{R}P_q^2) = \pi_{\widetilde{\omega}}(C(\mathbb{R}P_q^2))''$ may be identified with $B(l^2(\mathbb{N}))\cong d(B(l^2(\mathbb{N})))\subseteq C(\mathbb{Z}_2,B(l^2(\mathbb{N})))$, the space of constant functions from $\mathbb{Z}_2$ to $B(l^2(\mathbb{N}))$.\\

\noindent As it is well-known (and easy to show) that any ergodic action of $SU_q(2)$ on a C$^*$-algebra can be completed to a coaction of $\mathscr{L}^{\infty}(SU_q(2))$ on the von Neumann algebraic completion of the C$^*$-algebra in its GNS-representation with respect to the action-invariant state on it, we then obtain, from the ordinary algebraic definition of the actions of $SU_q(2)$ on the equatorial Podle\'{s} sphere (cf.~ \cite{Pod1},\cite{Mas1}) and on the quantum projective plane, the following von Neumann algebraic descriptions.

\begin{Def}\label{DefCoacProj} We denote by $\widetilde{\alpha}$ the unique coaction of $\mathscr{L}^{\infty}(SU_q(2))$ on $\mathscr{L}^{\infty}(S_{q\infty}^2)\cong C(\mathbb{Z}_2,B(l^2(\mathbb{N})))$ such that \[\widetilde{\alpha}(Y^*) = (a_+^*)^2\otimes Y^*-q(1+q^2)a_+^*b_+\otimes W - qb_+^2\otimes Y.\]\[\widetilde{\alpha}(W) = a_+^*b_+^*\otimes Y^* + (1-(1+q^2)b_+^*b_+)\otimes W + b_+a_+\otimes Y,\] \[\widetilde{\alpha}(Y) =  - q(b_+^*)^2\otimes Y^*-q(1+q^2)b_+^*a_+\otimes W+a_+^2\otimes Y ,\]This coaction is then $\mathbb{Z}_2$-equivariant.\\

\noindent We denote by $\alpha$ the restriction of $\widetilde{\alpha}$ to a coaction on $\mathscr{L}^{\infty}(\mathbb{R}P_q^2) = B(l^2(\mathbb{N}))\cong d(B(l^2(\mathbb{N})))$.
\end{Def}

\noindent Note then in particular that $\alpha$ is a coaction of the form treated in the final part of the previous section.

\subsection{Spectral decomposition of the $SU_q(2)$-action on the quantum projective plane}\label{SecPod}

\noindent In the previous section, we showed that the action of $SU_q(2)$ on the quantum projective plane $\mathbb{R}P_q^2$ gives rise to a coaction \[\alpha: B(l^2(\mathbb{N}))\rightarrow \mathscr{L}^{\infty}(SU_q(2))\bar{\otimes} B(l^2(\mathbb{N})),\] which is hence a projective corepresentation of $\mathscr{L}^{\infty}(SU_q(2))$ on $l^2(\mathbb{N})$ in the terminology of Definition \ref{DefProj}. In this section, we want to find an explicit description of the associated unitary $\mathcal{G}$ which we introduced in Definition \ref{DefProjCorep}.\\

\noindent Denote again by $e_{00}^{(+)}$ and $e_{00}^{(-)}$ the matrix units at position $00$ in resp.~ the $+$ and $-$ fiber of $B(l^2(\mathbb{N}))\oplus B(l^2(\mathbb{N}))$. Then $e_{00}^{(\pm)}$ are precisely the spectral projections at eigenvalue $\pm 1$ of $W$. Hence to determine $\alpha(e_{00})$, we should try to determine the eigenvectors for $1$ and $-1$ of $\widetilde{\alpha}_+(W)$, using the notation introduced at the end of the previous section (Notation \ref{Note}). To make the enunciation of the following Proposition and subsequent ones more succinct, the following notation, extending the one in Notation \ref{NotSUq11}, will come in handy.

\begin{Not} We denote $J^{(+)}=I_-^{(+)}$ for the set $\mathbb{Z}$. We denote $J^{(-)}=-I_{-}^{(-)}-1$ for the set $\mathbb{N}$. We denote by $J$ the disjoint union $J=J^{(+)}\sqcup J^{(-)}$.
\end{Not}

\noindent We then use the same notational conventions for elements in $J$ as for elements in $I_-$.\\

\begin{Prop}\label{PropEigvect} The spectrum of the operator $\widetilde{\alpha}_+(W)$ equals $\{\pm q^{2r},0\mid r\in \mathbb{N}\}$, with $0$ not occurring in the point spectrum.\\

\noindent For $r\in \mathbb{N}$, an orthonormal basis for the eigenspace of $\pm q^{2r}$ is given by the vectors $\xi_{r,\pm}^{(t,p)}$ with $p,t\in \mathbb{Z}$ and $p+r\in J^{(\pm)}$, determined by the formula \[ \xi_{r,\pm}^{(t,p)} = \sum_{n=0}^{\infty} Q_{q^{2}}(p_{\pm},r,n)\,e_n\otimes e_{t-n}\otimes e_{p+n} \in \mathscr{H}_+\otimes l^2(\mathbb{N}),\] with \begin{eqnarray*}Q_{q^2}(p_{\pm},r,n) &=& (\mp q)^r(\pm 1)^nq^{\frac{n(n-1)}{2}}  \frac{(\mp q^{2p+2r+2};q^2)_n(q^{4p+4n+4};q^4)_{\infty}^{1/2}}{(\mp q^{2p+2r+2};q^2)_{\infty}^{1/2}(q^4;q^4)_r^{1/2}(-1;q^2)_{\infty}^{1/2} (q^2;q^2)_n^{1/2}}\\ && \hspace{4.5cm} \cdot\qortb{q^{-2n}}{q^{-2r}}{\pm q^{-2p-2n}}{\mp q^{-2p-2n-2r}}{0}{q^2}{q^2}.\end{eqnarray*}

\end{Prop}

\noindent The proof of this Proposition will be presented in Appendix \ref{ApB}.\\

\noindent Now denote by $\mathscr{I}=\mathscr{I}_{\{-1,1\}}=\mathscr{I}_1+\mathscr{I}_{-1}$ the range of the spectral projection of $\widetilde{\alpha}_+(W)$ associated to the set $\{-1,1\}$. This then equals the range space of $\alpha(e_{00})$. Further recall that $I_-$ denotes the set $I_-^{(+)} \sqcup I_-^{(-)} = \mathbb{Z}\sqcup \mathbb{N}_0^-$. Then, by the above results, we can define a unitary map \[ u: \mathscr{I} \rightarrow \mathscr{H}_-=l^2(I_-)\otimes l^2(\mathbb{Z}):\xi_{0,\pm}^{(t,p)} \rightarrow e_{-p-1}^{(\pm)}\otimes e_{p+t},\] where we recall that $e_n^{(\pm)}=e_{n_{\pm}}$ for $\mathbf{n}\in I_-$.\\

\noindent In the following Proposition, we will again use the notation $I_+ = \mathbb{N}$.

\begin{Prop}\label{PropIdenG} The map \[\mathcal{G}:\mathscr{H}_-\otimes l^2(\mathbb{N})\rightarrow \otimes \mathscr{H}_-\otimes l^2(\mathbb{N}):\xi_{r,\pm}^{(t,p)} \rightarrow e_{-p-r-1}^{(\pm)}\otimes e_{p+t-r}\otimes e_r\] defines a unitary, and $(u^*\otimes 1)\mathcal{G}$ coincides with the unitary constructed in Definition \ref{DefProjCorep}.\end{Prop}

\begin{proof} The fact that $\mathcal{G}$ is a well-defined unitary is of course immediate by the previous Proposition. Then, from the way $\widetilde{\alpha}_+(Y^*)$ acts on the non-normalized eigenvectors $\eta_{r,\pm}^{(t,p)}$ in the proof of Proposition \ref{PropEigvect} (see the identities (\ref{EqActY}), (\ref{EqActY2}) in Appendix \ref{ApB}), we have, for $r\in\mathbb{N}$, $t\in \mathbb{Z}$ and $p\in J^{(\pm)}$, that \begin{eqnarray*} \mathcal{G}^*(u\otimes 1) \, \xi_{0,\pm}^{(t,p)}\otimes e_r &=& \mathcal{G}^*\, e_{-p-1}^{(\pm)} \otimes e_{p+t}\otimes e_r \\ &=& \xi_{r,\pm}^{(t+2r,p-r)}\\ &=& q^{-r(r+1)}(q^4;q^4)_r^{-1/2} \widetilde{\alpha}_+((YW)^*)^r \xi_{0,\pm}^{(t,p)} \\ &=& q^{-r(r+1)}(q^4;q^4)_r^{-1/2} \alpha((YW)^*)^r \xi_{0,\pm}^{(t,p)} \\ &=& \alpha(e_{r0}) \xi_{0,\pm}^{(t,p)},\end{eqnarray*} which proves the Proposition.\end{proof}

\noindent As mentioned in the first section, we may treat $\mathcal{G}$ itself as the implementing unitary of $\alpha$, as the unitary $u$ only serves to reparametrize the Hilbert space $\alpha(e_{00})(\mathscr{H}_+\otimes l^2(\mathbb{N}))$.\\

\noindent Now further denote \[e:l^2(I_-)\otimes l^2(\mathbb{Z}) \rightarrow l^2(I_-)\otimes l^2(\mathbb{Z}): e_{n}^{(\pm)} \otimes e_{l}\rightarrow \pm \, e_{n}^{(\pm)}\otimes e_l,\] and denote $P_{\pm} = \frac{1}{2}(1\pm e)$ and \[\mathcal{G}^{(\pm)} = (P_{\pm}\otimes 1)\mathcal{G}.\] Then $e$ is precisely the self-adjoint unitary which appeared in Proposition \ref{PropGL}, so our notation is consistent with the one introduced in Notation \ref{NotRest}. It is then also easy to see that \[\mathcal{G}^{(\pm)}\widetilde{\alpha}_+(W) = \pm (1\otimes W^{(+)})\mathcal{G}^{(\pm)},\]\[ \mathcal{G}^{(\pm)}\widetilde{\alpha}_+(Y) = \pm (1\otimes Y^{(+)})\mathcal{G}^{(\pm)}.\] Hence \[\mathcal{G}^{(\pm)}\widetilde{\alpha}_+(x) = \pm (1\otimes x^{(+)})\mathcal{G}^{(\pm)}, \qquad \textrm{for all }x\in \mathscr{L}^{\infty}(S_{q\infty}^2).\]

\noindent Write now \[\mathcal{G}^{(\pm)} = \sum_{r,s=0}^{\infty} \mathcal{G}_{r,s}^{(\pm)}\otimes e_{rs}\] as a $\sigma$-weakly converging sum, where $\mathcal{G}_{r,s}^{(\pm)}: l^2(I_+)\otimes l^2(\mathbb{Z})\rightarrow l^2(I_-^{(\pm)})\otimes l^2(\mathbb{Z})$.

\begin{Prop}\label{PropFormG} For $n\in \mathbb{N},k\in \mathbb{Z}$, we have \[\!\!\!\!\!\!\!\!\!\!\!\!\!\mathcal{G}_{r,s}^{(\pm)} e_{n}\otimes e_k = \frac{1}{\sqrt{2}}(\mp q)^r (\pm 1)^n q^{\frac{n(n-1)}{2}} \frac{(\mp q^{2s+2r-2n+2};q^2)_{n}^{1/2}(q^{2n+2};q^2)_{\infty}^{1/2}}{(\mp q^{2s+2r+2};q^2)_{\infty}^{1/2} (q^4;q^4)_r^{1/2}(q^4;q^4)_s^{1/2}} \]\[\qquad \qquad \qquad \cdot \qortb{q^{-2n}}{q^{-2r}}{\pm q^{-2s}}{\mp q^{-2s-2r}}{0}{q^2}{q^2}\, e_{n-r-s-1}^{(\pm)}\otimes e_{k-r+s}.\]
\end{Prop}

\noindent \emph{Remark:} Using transformation formula (III.11) of \cite{Gas1}, we have \[\qortb{q^{-2n}}{q^{-2r}}{-q^{-2s}}{q^{-2s-2r}}{0}{q^2}{q^2} = (-1)^n \qortb{q^{-2n}}{-q^{-2r}}{q^{-2s}}{q^{-2s-2r}}{0}{q^2}{q^2}.\] Hence we can also write  \[\!\!\!\!\!\!\!\!\!\!\!\!\!\mathcal{G}_{r,s}^{(\pm)} e_{n}\otimes e_k = \frac{1}{\sqrt{2}}(\mp q)^r q^{\frac{n(n-1)}{2}} \frac{(\mp q^{2s+2r-2n+2};q^2)_{n}^{1/2}(q^{2n+2};q^2)_{\infty}^{1/2}}{(\mp q^{2s+2r+2};q^2)_{\infty}^{1/2} (q^4;q^4)_r^{1/2}(q^4;q^4)_s^{1/2}} \]\[\qquad \qquad \qquad \cdot \qortb{q^{-2n}}{\pm q^{-2r}}{q^{-2s}}{\mp q^{-2s-2r}}{0}{q^2}{q^2}\, e_{n-r-s-1}^{(\pm)}\otimes e_{k-r+s}.\]

\begin{proof} We have that for $n\in \mathbb{N}$, $k,l\in \mathbb{Z}$ and $m\in I_-^{(\pm)}$, \begin{eqnarray*} \langle e_m^{(\pm)}\otimes e_l,\mathcal{G}_{r,s}^{(\pm)}\, e_n\otimes e_k\rangle &=& \langle e_m^{(\pm)}\otimes e_l\otimes e_r,\mathcal{G}^{(\pm)}\, e_n\otimes e_k\otimes e_s\rangle \\ &=& \langle (\mathcal{G}^{(\pm)})^*\,e_m^{(\pm)}\otimes e_l\otimes e_r,e_n\otimes e_k\otimes e_s\rangle\\ &=& \langle  \xi_{r,\pm}^{(l+2r+m+1,-m-r-1)},e_n\otimes e_k\otimes e_s\rangle \\ &=& \delta_{m,n-r-s-1}\delta_{l,k-r+s}\langle  \xi_{r,\pm}^{(k+n,s-n)},e_n\otimes e_k\otimes e_s\rangle \end{eqnarray*} The Proposition then follows immediately by the concrete form of the $\xi_{r,\pm}^{(t,p)}$ given in Proposition \ref{PropEigvect}.

\end{proof}

\begin{Def}\label{DefOp} We define the following operators: \[\begin{array}{lll} L_{0+}:& l^2(I_+)\otimes l^2(\mathbb{Z})\rightarrow l^2(I_0)\otimes l^2(\mathbb{Z}):&e_{n}\otimes e_k \rightarrow (q^{2n+2};q^2)_{\infty}^{1/2} e_n\otimes e_k,\\ a_0:&l^2(I_0)\otimes l^2(\mathbb{Z}) \rightarrow l^2(I_0)\otimes l^2(\mathbb{Z}): &e_n\otimes e_k \rightarrow e_{n-1}\otimes e_k,\\ L_{-0}^{(\pm)}:& l^2(I_0)\otimes l^2(\mathbb{Z})\rightarrow  l^2(I_-)\otimes l^2(\mathbb{Z}):&e_n\otimes e_k \rightarrow  q^{\frac{n(n+1)}{2}} (\mp q^{-2n};q^2)_{\infty}^{1/2} e_n^{(\pm)}\otimes e_k, \\ f:& l^2(I_-)\otimes l^2(\mathbb{Z})\rightarrow  l^2(I_-)\otimes l^2(\mathbb{Z}): &e_{n}^{(\pm)}\otimes e_k \rightarrow (\pm 1)^{n+1} e_{n}^{(\pm)}\otimes e_k. \end{array}\]

\end{Def}

\noindent Note that the map $L_{0+}$ also appeared in \cite{DeC5} (and \cite{Wor6}).

\begin{Def} For $r,s\in \mathbb{N}$, we define polynomials $K_{r,s}^{(\pm)}$ as follows: \[K_{r,s}^{(\pm)}(x) = \,\!_2\varphi_1\left(\begin{array}{cc} q^{-2\min\{r,s\}} & -q^{-2\min\{r,s\}} \\ &\!\!\!\!\!\!\!\!\!\!\!\!\!\!\!\!\!\!\!\!\!\!\!\!\!\!\!\pm q^{2|r-s|+2} \end{array}\mid q^2,q^2x\right).\]

\end{Def}

\begin{Prop} For $s\geq r$, we have
\[\mathcal{G}_{r,s}^{(\pm)} = \frac{(\pm 1)^{s}}{\sqrt{2}} q^{\frac{1}{2}(r-s)(3r+s+1)} \frac{(q^4;q^4)_s^{1/2}}{(q^4;q^4)_r^{1/2}} \frac{(\pm q^{2s-2r+2};q^2)_{\infty}}{(q^4;q^4)_{\infty}} L_{-0}^{(\pm)} a_0^{r+s+1} fL_{0+} K_{r,s}^{(\pm)}(b^*b) b^{s-r}.\]

\noindent For $r\geq s$, we have \[\mathcal{G}_{r,s}^{(\pm)} = \frac{(\pm 1)^{r}}{\sqrt{2}} q^{\frac{1}{2}(s-r)(3s+r+1)} \frac{(q^4;q^4)_r^{1/2}}{(q^4;q^4)_s^{1/2}} \frac{(\pm q^{2r-2s+2};q^2)_{\infty}}{(q^4;q^4)_{\infty}} L_{-0}^{(\pm)} a_0^{r+s+1} L_{0+} K_{r,s}^{(\pm)}(b^*b) (- q b^*)^{r-s}.\]
\end{Prop}

\begin{proof}

\noindent For $s\geq r$, we have, by applying the transformation formula (III.6) of \cite{Gas1} with respect to $q^{-2r}$ as the terminating factor, that \[\qortb{q^{-2n}}{q^{-2r}}{\pm q^{-2s}}{\mp q^{-2s-2r}}{0}{q^2}{q^2} = \]\[ (-q^{2r-2n})^r \frac{(\pm q^{2s-2r+2};q^2)_r}{(\mp q^{2s+2};q^2)_r} \qorta{q^{-2r}}{-q^{-2r}}{\pm q^{2s-2r+2}}{q^2}{q^{2n+2}},\] while for $r\geq s$, we have, applying transformation formula (III.6) of \cite{Gas1}  with respect to $q^{-2s}$ as the terminating factor, \[\qortb{q^{-2n}}{\pm q^{-2r}}{q^{-2s}}{\mp q^{-2s-2r}}{0}{q^2}{q^2} =\]
\[ (-q^{2s-2n})^s \frac{(\pm q^{2r-2s+2};q^2)_s}{(\mp q^{2r+2};q^2)_s} \qorta{q^{-2s}}{-q^{-2s}}{\pm q^{2r-2s+2}}{q^2}{q^{2n+2}}.\]
Then with the above transformation formulas at hand, the Proposition follows straightforwardly from the formulas for $\mathcal{G}_{r,s}^{(\pm)}$ in Proposition \ref{PropFormG} and the remark following it.
\end{proof}

\noindent \emph{Remark:} It seems odd that in the formula for $\mathcal{G}_{r,s}^{(\pm)}$, an extra `parity operator' $f$ appears when switching from $r\geq s$ to $s\geq r$. We have no real conceptual reason to explain this phenomenon.

\section{Identification of the reflected quantum group}

\noindent In the previous section, we found an explicit description of the unitary $\mathcal{G}$ implementing the projective corepresentation $\alpha$ of $\mathscr{L}^{\infty}(SU_q(2))$ on $B(l^2(\mathbb{N}))\cong \mathscr{L}^{\infty}(\mathbb{R}P_q^2)$. As $\mathcal{G}$ was defined as a unitary from $\mathscr{H}_+=l^2(I_+)\otimes l^2(\mathbb{Z})$ to $\mathscr{H}_-=l^2(I_-)\otimes l^2(\mathbb{Z})$, we will have that the space $N$ for $\mathcal{G}$, introduced in the Notation \ref{NotCoob}, will be a subspace of $B(\mathscr{H}_+,\mathscr{H}_-)$.

\begin{Lem} The equality $N=\mathscr{L}(-,+)$ holds.\end{Lem}

\begin{proof} We recall that $\mathscr{L}(-,+)$ was just the space $B(l^2(I_+),l^2(I_-))\bar{\otimes}\mathscr{L}(\mathbb{Z}) \subseteq B(\mathscr{H}_+,\mathscr{H}_-)$. As it is immediately observed from Proposition \ref{PropFormG} that all $\mathcal{G}_{0,r}=\mathcal{G}_{0,r}^{(+)}+\mathcal{G}_{0,r}^{(-)}$ commute with $1\otimes S$, with $S$ the bilateral forward shift on $l^2(\mathbb{Z})$, it follows that $\mathcal{G}_{0,r}\subseteq \mathscr{L}(-,+)$ for all $r\in \mathbb{N}$. As $\mathscr{L}^{\infty}(SU_q(2))=\mathscr{L}(+,+)=B(l^2(I_+))\bar{\otimes}\mathscr{L}(\mathbb{Z})$ by definition, and $N$ is generated by the $\mathcal{G}_{0,r}$ as a right $\mathscr{L}(+,+)$-module (again by definition), we obtain the inclusion $N\subseteq \mathscr{L}(-,+)$. \\

\noindent Now, using the notation introduced at the end of the previous section, we have $\mathcal{G}_{0,r}^{(\pm)}\in N$ by Lemma \ref{LemPm}. Take $p\in \mathbb{N}$. It follows then from Proposition \ref{PropFormG}, that $\mathcal{G}_{0,0}^{(+)}(e_{p0}\otimes 1)$ is a non-zero scalar multiple of the matrix unit $e_{(p-1)_+,0}\otimes 1$ in $B(\mathscr{H}_{+},\mathscr{H}_-)$, while $\mathcal{G}_{0,p}^{(\pm)}(e_{00}\otimes S^{-p})$ is a non-zero scalar multiple of the matrix unit $e_{(-p-1)_{\pm},0}\otimes 1$. So $N$ contains all matrix units of the form $e_{\mathbf{p},0}\otimes 1$, for $\mathbf{p}\in I$. As $N$ is by definition closed under right multiplication with elements in $\mathscr{L}(+,+)$, it follows that indeed $N=\mathscr{L}(-,+)$.
\end{proof}

\noindent It then follows immediately from this Lemma that the linking von Neumann algebra $\left(\begin{array}{ll} P & N \\ O & M \end{array}\right)$ associated to $\mathcal{G}$, using again Notation \ref{NotCoob}, equals precisely the $\{-,+\}$-part of the von Neumann algebra $Q$ in Theorem \ref{TheoCoa}, i.e.~ $\left(\begin{array}{ll} P & N \\ O & M \end{array}\right)=\left(\begin{array}{ll} \mathscr{L}(-,-) & \mathscr{L}(-,+) \\ \mathscr{L}(+,-) & \mathscr{L}(+,+) \end{array}\right)$, with in particular $M=\mathscr{L}^{\infty}(SU_q(2))$. We also recall that we denoted by $\Gamma$ and $\Gamma_{\mu\nu}$ (or $\Gamma_M,\Gamma_N,\ldots$) the comultiplication and its constituents on $Q$, as obtained by the method explained in Proposition \ref{TheoTwist}. Further recall that we defined another collection of maps $\Delta_{\mu\nu}$ on $\mathscr{L}(\mu,\nu)$ in the discussion preceding Theorem \ref{TheoCoa}.

\begin{Prop}\label{PropEqCom} The comultiplications $\Delta_{\mu\nu}$ and $\Gamma_{\mu\nu}$ coincide for $\mu,\nu\in\{-,+\}$.

\end{Prop}

\noindent From this Proposition (and from Proposition \ref{TheoTwist}), it will then immediately follow that the $\Delta_{\mu\nu}$ for $\mu,\nu\in \{-,+\}$ are coassociative. In particular, as we will not actually need to use that $\Delta_{-}$ is coassociative, this will give an alternative proof for the coassociativity of the comultiplication $\Delta_{-}$ on $\mathscr{L}^{\infty}(\widetilde{SU}_q(1,1))$ (first established in \cite{Koe3}).\\

\begin{proof}[Proof (of Proposition \ref{PropEqCom})] The fact that $\Gamma_{++}=\Delta_{++}$, the natural comultiplication on $\mathscr{L}^{\infty}(SU_q(2))$, was part of Proposition \ref{TheoTwist}.\\

\noindent Further, as $\mathscr{L}(+,-) = \mathscr{L}(-,+)^*$, and $\mathscr{L}(-,+)\cdot \mathscr{L}(+,-)$ is $\sigma$-weakly dense in $\mathscr{L}(-,-)$, the equalities $\Delta_{--}=\Gamma_{--}$ and $\Delta_{+-}=\Gamma_{+-}$ will immediately follow from the equality $\Delta_{-+}=\Gamma_{-+}$.\\

\noindent But $\Delta_{-+}(xy) = \Delta_{-+}(x)\Delta_{++}(y)$ and $\Gamma_{-+}(xy)=\Gamma_{-+}(x)\Gamma_{++}(y)$ for all elements $x\in \mathscr{L}(-,+)$ and $y\in \mathscr{L}(+,+)$. So we see that it is sufficient to prove, for $\mathbf{p}\in I_-$, that $\Delta_N(e_{\mathbf{p}0}\otimes 1) = \Gamma_{N}(e_{\mathbf{p}0}\otimes 1)$. As both these expressions vanish on $\xi_{r,s,p,t}^+$ for $p\in \mathbb{N}_0$, it is already sufficient to prove that, for $\mathbf{p}\in I_-= \mathbb{Z}\sqcup \mathbb{N}_0^-$ and $r,s,t\in \mathbb{Z}$, we have \[\Gamma_{N}(e_{\mathbf{p}0}\otimes 1) \xi_{r,s,0,t}^+ = \xi_{r,s,\mathbf{p},t}^{-}.\]

\noindent Using the formulas for the $\mathcal{G}_{r,s}$ in Proposition \ref{PropFormG}, we see that for $p\in \mathbb{N}$, \[\mathcal{G}_{0,0}^{(+)}(e_{p0}\otimes 1) = \frac{1}{\sqrt{2}} q^{\frac{1}{2}p(p-1)} \frac{(-q^{-2p+2};q^2)_p^{1/2}(q^{2p+2};q^2)_{\infty}^{1/2}}{(-q^2;q^2)_{\infty}^{1/2}} \,e_{(p-1)_+,0}\otimes 1\] and \[ \mathcal{G}_{0,p}^{(\pm)}(e_{00}\otimes 1) = \frac{1}{\sqrt{2}} \frac{(q^2;q^2)_{\infty}^{1/2}}{(\mp q^{2p+2};q^2)_{\infty}^{1/2}(q^4;q^4)_{p}^{1/2}} \, e_{(-p-1)_{\pm},0}\otimes S^{p}.\]

\noindent Then since, for $m,n,p\in \mathbb{N}$ and $r,s,t,k\in \mathbb{Z}$, we have $\Gamma_{++}(e_{mn}\otimes S^k)\xi_{r,s,p,t}^+ = \delta_{n,p}\, \xi_{r,s,m,t+k}^+$ by Proposition \ref{PropBasSUq2}, and $\Gamma_N(xy)= \Gamma_N(x)\Gamma_M(y)$ for $x\in N$ and $y\in M=\mathscr{L}^{\infty}(SU_q(2))$, we see that it is sufficient to prove that for $p\in \mathbb{N}_0$, we have \begin{equation}\label{eqn1} \xi_{r,s,(p-1)_+,t}^- = \sqrt{2} q^{-\frac{1}{2}p(p-1)} \frac{(-q^2;q^2)_{\infty}^{1/2}}{(-q^{-2p+2};q^2)_p^{1/2}(q^{2p+2};q^2)_{\infty}^{1/2}} \, \Gamma_{N}(\mathcal{G}_{0,0}^{(+)})\xi_{r,s,p,t}^+,\end{equation} and that for $p\in \mathbb{N}$, we have \begin{equation}\label{eqn2}\xi_{r,s,(-p-1)_{\pm},t+p}^{-} =  \sqrt{2} \frac{(\mp q^{2p+2};q^2)_{\infty}^{1/2}(q^4;q^4)_{p}^{1/2}}{(q^2;q^2)_{\infty}^{1/2}}\,\Gamma_{N}(\mathcal{G}_{0,p}^{(\pm)})\xi_{r,s,0,t}^+.\end{equation} \\

\noindent From Lemma \ref{LemPm}, it follows that \[\Gamma_N(\mathcal{G}_{r,s}^{(+)}) = \sum_{\pm}\sum_{j\in \mathbb{N}} \mathcal{G}_{r,j}^{(\pm)}\otimes \mathcal{G}_{j,s}^{(\pm)}\] and  \[\Gamma_N(\mathcal{G}_{r,s}^{(-)}) = \sum_{\pm}\sum_{j\in \mathbb{N}} \mathcal{G}_{r,j}^{(\mp)}\otimes \mathcal{G}_{j,s}^{(\pm)}\] as $\sigma$-weakly converging sums.\\

\noindent We will in the following use the alternative expression (\ref{eqnPfunc}) for the function $P_{q^2}^+$ which appears in Proposition \ref{PropBasSUq2}. Then using the formulas from Proposition \ref{PropFormG}, we get from the formulas in the previous paragraph, after some easy simplifications, \begin{eqnarray*}&&\Gamma_N(\mathcal{G}_{0,0}^{(+)})\xi_{r,s,p,t}^+\;=\\ &&\qquad \frac{1}{2} \sum_{\pm} \underset{v-w=t}{\sum_{v,w\in I_+}}\sum_{j\in \mathbb{Z}} (-1)^p (\pm 1)^{v} (\mp q)^j  q^{vw+p(v+w+1)} q^{\frac{1}{2}v(v-1)}q^{\frac{1}{2}w(w-1)} \\ &&\qquad \cdot \frac{(q^{2v+2};q^2)_{\infty} (q^{2p+2};q^2)_{\infty}^{1/2}(\mp q^{2j-2v+2};q^2)_{v}^{1/2}(\mp q^{2j-2w+2};q^2)_w^{1/2}(\pm q^{2j+2};q^2)_{\infty}}{(q^2;q^2)_w(q^4;q^4)_{\infty}}\\ &&\qquad \cdot \qortb{q^{-2w}}{q^{-2v}}{q^{-2p}}{0}{0}{q^2}{q^2}\, e_{v-j-1}^{(\pm)} \otimes e_{r+p-w+j} \otimes e_{w-j-1}^{(\pm)}\otimes e_{s+v-p-j},\end{eqnarray*} where we recall that by convention $e_{\mathbf{a}}=0$ when $\mathbf{a}\notin I_-$. Now we apply the change of variables $j=w+m$. Changing then the order of summation, we get, using still the notation $J^{(+)} = \mathbb{Z}$ and $J^{(-)} = \mathbb{N}$, \begin{eqnarray*}&&\Gamma_N(\mathcal{G}_{0,0}^{(+)})\xi_{r,s,p,t}^+\;=\\ &&\qquad \frac{1}{2} (-1)^p q^{\frac{1}{2}(2tp+2p+t^2-t)} \frac{(q^{2p+2};q^2)_{\infty}^{1/2}}{(q^4;q^4)_{\infty}} \sum_{\pm} (\pm 1)^t  \underset{m-n=t}{\sum_{m,n\in J^{(\pm)}}} (\mp q)^m \underset{v-w=t}{\sum_{v,w\geq 0}} (-1)^w q^{w(2v+2p)}\\ && \qquad \cdot \frac{(q^{2v+2};q^2)_{\infty} (\mp q^{2n+2};q^2)_{v}^{1/2} (\mp q^{2m+2};q^2)_w^{1/2}(\pm q^{2m+2w+2};q^2)_{\infty}}{(q^2;q^2)_w}\\ &&\qquad \cdot \qortb{q^{-2w}}{q^{-2v}}{q^{-2p}}{0}{0}{q^2}{q^2} e_{-n-1}^{(\pm)} \otimes e_{r+p+m}\otimes e_{-m-1}^{(\pm)}\otimes e_{s-n-p}.\end{eqnarray*}

\noindent Now with the conditions on $m,n,v,w$ as in this summation, we can write \[ (\mp q^{2n+2};q^2)_v^{1/2} (\mp q^{2m+2};q^2)_w^{1/2} = \frac{(\mp q^{2n+2};q^2)_{\infty}^{1/2}}{(\mp q^{2m+2};q^2)_{\infty}^{1/2}}(\mp q^{2m+2};q^2)_w.\] Then we can apply the summation formula in Proposition \ref{PropSum1} with $x=\pm q^{2m}$ and $y=\pm q^{2n}$ to obtain \begin{eqnarray*} &&\Gamma_N(\mathcal{G}_{0,0}^{(+)})\xi_{r,s,p,t}^+ \,=\\ &&\qquad \frac{1}{2} (-1)^p q^{\frac{1}{2}(2tp+2p+t^2-t)}\frac{(-1;q^2)_{p}(q^{2p+2};q^2)_{\infty}^{1/2}}{(q^4;q^4)_{\infty}}\sum_{\pm} (\pm 1)^t\underset{m-n=t}{\sum_{m,n\in J^{(\pm)}}} (\mp q)^m\\&&\qquad  \cdot \frac{(\mp q^{2n+2};q^2)_{\infty}^{1/2}}{(\mp q^{2m+2};q^2)_{\infty}^{1/2}} \qortPsi{\mp q^{-2n}}{q^{2t+2}}{q^2}{\pm q^{2m+2p+2}} e_{-n-1}^{(\pm)}\otimes e_{r+p+m}\otimes e_{-m-1}^{(\pm)} \otimes e_{s-n-p}.\end{eqnarray*} Plugging this into the right hand side of equality (\ref{eqn1}), a straightforward comparison of coefficients proves the validity of the identity (\ref{eqn1}).\\

\noindent The other case $(\ref{eqn2})$ is entirely similar. We first write $P_{q^{2}}^+$ again in the $_3\varphi_2$-form as in the previous case, as this then simplifies immediately inside the expression for $\xi_{r,s,0,t}^+$. Next, for $\mu\in \{-,+\}$, write the expansion for $\Gamma_N(\mathcal{G}_{0,p}^{(\mu)})$ as \[\Gamma_N(\mathcal{G}_{0,p}^{(\mu)}) = \sum_{\pm}\sum_{j\in \mathbb{N}} \mathcal{G}_{0,j}^{(\pm \mu)}\otimes \mathcal{G}_{j,p}^{(\pm)},\] and, when evaluating it on $\xi_{r,s,0,t}^+$, use the formula in Proposition \ref{PropFormG} for the factor $\mathcal{G}_{0,j}^{(\pm\mu)}$, and the slightly different expression given in the remark below that Proposition for the factor $\mathcal{G}_{j,p}^{(\pm)}$. We find \begin{eqnarray*} &&\Gamma_N(\mathcal{G}_{0,p}^{(\mu)}) \xi_{r,s,0,t}^{+} =\\&&\qquad  \frac{1}{2} \frac{1}{(q^4;q^4)_{\infty}(q^2;q^2)_{\infty}^{1/2}(q^{4};q^4)_p^{1/2}}\sum_{\pm} \underset{v-w=t}{\sum_{v,w\in I^+}}\;\underset{i-j=p}{\sum_{i,j\in \mathbb{Z}}} (\pm\mu)^v(\mp 1)^j q^{vw+\frac{1}{2}v(v-1)+\frac{1}{2}w(w-1)+j}\\ &&\qquad \frac{(\mp q^{2i-2w+2};q^2)_{\infty}^{1/2}(\mp \mu q^{2j-2v+2};q^2)_{\infty}^{1/2}(q^{2v+2};q^2)_{\infty}(q^{2w+2};q^2)_{\infty} (\pm \mu q^{2j+2};q^2)_{\infty}} {(\mp q^{2i+2};q^2)_{\infty}} \\&& \qquad \qortb{q^{-2w}}{\pm q^{-2j}}{q^{-2p}}{\mp q^{-2i}}{0}{q^2}{q^2} e_{v-j-1}^{(\pm \mu)}\otimes e_{r-w+j}\otimes e_{w-i-1}^{(\pm)}\otimes e_{s+v-j+p}.\end{eqnarray*} If we then apply the change of variables $j=v+m$ and $i=w+n$ and change the order of summation, we obtain \begin{eqnarray*} &&\Gamma_N(\mathcal{G}_{0,p}^{(\mu)}) \xi_{r,s,0,t}^{+} =\\&&\qquad  \frac{1}{2} \frac{(-\mu)^tq^{\frac{1}{2}t(t+1)}}{(q^4;q^4)_{\infty}(q^2;q^2)_{\infty}(q^{4};q^4)_p^{1/2}}\sum_{\pm} \underset{m-n=-t-p}{\sum_{m,n\in J^{(\pm)}}} (\mp q)^m (\mp \mu q^{2m+2},\mp q^{2n+2};q^2)_{\infty}^{1/2} \\ && \qquad \sum_{w=0}^{\infty} (-\mu)^w q^{2w^2} q^{(2n-2m-2p)w}\frac{(\pm \mu q^{2n-2p+2w+2},q^{2n-2m-2p+2w+2},q^{2w+2};q^2)_{\infty}}{(\mp q^{2n+2w+2};q^2)_{\infty}}\\ && \qquad \qortb{q^{-2w}}{\pm q^{-2n+2p-2w}}{q^{-2p}}{\mp q^{-2n-2w}}{0}{q^2}{q^2} e^{(\pm \mu)}_{-m-1} \otimes e_{r+n-p} \otimes e_{-n-1}^{(\pm)}\otimes e_{s-m+p}.\end{eqnarray*}

\noindent Using then the summation formula in Proposition \ref{PropSum2} with $x=\pm q^{2n-2p}$ and $y= \pm q^{2m}$ (and with the $\pm$ in that Proposition replaced by $\mu$), and plugging this expression for $\Gamma_N(\mathcal{G}_{0,p}^{(\mu)}) \xi_{r,s,0,t}^{+}$ into the right hand side of the identity (\ref{eqn2}), we can again conclude that both sides are equal.

\end{proof}

\section{The linking weak von Neumann bialgebra between $SU_q(2)$ and $\widetilde{E}_q(2)$}

\noindent The main part of this section consists in showing that, in Theorem \ref{TheoCoa}, the $\Delta_{\mu\nu}$ with $\mu,\nu\in \{+,0\}$ are coassociative. This will then let us conclude the proof of that Theorem in a straightforward manner.\\

\noindent For the $\{0,+\}$-part of $(Q,\Delta_Q)$, we do not have to go through as much trouble as for the $\widetilde{SU}_q(1,1)$-case, as we have already treated a lot of material concerning it in \cite{DeC5}. Let us state one of the main results from that paper (see Theorem 3.13 and Theorem 4.4). We will use notations as in the introduction. Also recall that the operator $L_{0+}$ we will use was defined in Definition \ref{DefOp}.

\begin{Prop}  There exists a unital, normal, coassociative $^*$-homomorphism \[\Gamma: \left(\begin{array}{cc} \mathscr{L}(0,0) & \mathscr{L}(0,+)\\ \mathscr{L}(+,0)& \mathscr{L}(+,+)\end{array}\right) \rightarrow \left(\begin{array}{cc} \mathscr{L}(0,0)\bar{\otimes} \mathscr{L}(0,0) & \mathscr{L}(0,+)\bar{\otimes}\mathscr{L}(0,+)\\ \mathscr{L}(+,0)\bar{\otimes} \mathscr{L}(+,0)& \mathscr{L}(+,+)\bar{\otimes}\mathscr{L}(+,+)\end{array}\right)\] which restricts to maps \[\Gamma_{\mu\nu}: \mathscr{L}(\mu,\nu)\rightarrow\mathscr{L}(\mu,\nu)\bar{\otimes}\mathscr{L}(\mu,\nu),\] and such that $\Gamma_{00} = \Delta_{\mathscr{L}^{\infty}(\widetilde{E}_q(2))}$, $\Gamma_{++} = \Delta_{\mathscr{L}^{\infty}(SU_q(2))}$ and \[\Delta_{0+}(L_{0+}) = \sum_{k=0}^{\infty} (q^2;q^2)_k^{-1} \;v_0^kL_{0+}b_+^k\otimes v_0^kL_{0+}(-qb_+^*)^k,\] this last sum being convergent in norm.
\end{Prop}

\noindent The coassociativity of the $\Delta_{\mu\nu}$ will then follow immediately from the following Proposition.

\begin{Prop}\label{PropIdLink} For all $\mu,\nu\in \{0,+\}$, we have $\Delta_{\mu\nu} = \Gamma_{\mu\nu}$.\end{Prop}

\begin{proof} We first remark that we have not substantiated yet the claim that $\Delta_0$, in the form represented in Proposition \ref{PropBasEq2}, equals the comultiplication in Definition \ref{DefEq2}, which we will momentarily write as $\Delta_{\mathscr{L}^{\infty}(\widetilde{E}_q(2))}$ for distinction. This equality will be proven in the course of the proof.\\

\noindent We observe that the Proposition will follow once we prove that $\Delta_0(a_0) = \Gamma_{00}(a_0)$ and $\Delta_{0+}(L_{0+})=\Gamma_{0+}(L_{0+})$. Indeed, suppose this is satisfied. Then as the equality \[(\Delta_{\mathscr{L}^{\infty}(SU_q(2))}=)\Delta_+=\Gamma_+\] was part of the previous Proposition, we will have that, for $x\in \mathscr{L}^{\infty}(SU_q(2))$ and $m\in \mathbb{Z}$, \begin{eqnarray*} \Delta_{0+}(a_0^mL_{0+}x) &=& \Delta_{0}(a_0)^m\Delta_{0+}(L_{0+})\Delta_+(x) \\ &=& \Gamma_{0}(a_0)^m \Gamma_{0+}(L_{0+})\Gamma_+(x) \\ &=& \Gamma_{0+}(a_0^mL_{0+}x).\end{eqnarray*} But the elements of the form $a_0^mL_{0+}x$ are easily seen to be $\sigma$-weakly dense in $\mathscr{L}(0,+)$. So then $\Delta_{0+}=\Gamma_{0+}$ follows from this. As in the beginning of Proposition \ref{PropEqCom}, this will allow us to conclude $\Delta_{\mu\nu}=\Gamma_{\mu\nu}$ for all $\mu,\nu\in \{0,+\}$.\\

\noindent We now first show that $\Delta_0(a_0)=a_0\otimes a_0$. By the definition of $\Delta_0$ as given in Proposition \ref{PropBasEq2}, this means that, for $r,s,t\in\mathbb{Z}$ and $p\in I_0=\mathbb{Z}$, we should have \[\xi_{r,s,p-1,t}^{0} = \underset{v-w=t}{\sum_{v,w\in I_0}} P_{q^2}^0(p,v,w) e_{v-1} \otimes e_{r+p-w}\otimes e_{w-1}\otimes e_{s-p+v}.\] Comparing coefficients, this reduces to the identity \[P_{q^{2}}^0(p-1,v-1,w-1) = P_{q^2}^0(p,v,w)\] for all $p,v,w\in I_0$. This is in fact immediate from the definition of $P_{q^2}^0$, as the formula only involves pairwise differences of the $v,w$ and $p$.\\

\noindent Thus the only point left to prove is that $\Delta_{0+}(L_{0+})=\Gamma_{0+}(L_{0+})$. Choose a vector $\xi_{r,s,p,t}^+$ with $r,s,t\in\mathbb{Z}$ and $p\in I_+=\mathbb{N}$. It is then sufficient to prove that $\Delta_{0+}(L_{0+})\xi_{r,s,p,t}^+ = \Gamma_{0+}(L_{0+})\xi_{r,s,p,t}^+$. Using the definition of $\Delta_{0+}$, and the formulas for $L_{0+}$ and $\Gamma_{0+}(L_{0+})$ stated in the previous Proposition, this becomes the identity \begin{eqnarray*} (q^{2p+2};q^2)_{\infty}^{1/2}\xi_{r,s,p,t}^0 &=& \underset{v-w=t}{\sum_{v,w\in I_+}}\sum_{k=0}^{\infty} P_{q^2}^+(p,v,w) (-q)^kq^{k(v+w)} \frac{(q^{2v+2},q^{2w+2};q^2)_{\infty}^{1/2}}{(q^2;q^2)_k} \\ && \qquad \qquad\qquad\cdot e_{v-k}\otimes e_{r+p-w+k}\otimes e_{w-k}\otimes e_{s-p+v-k}.\end{eqnarray*} Multiplying both sides with $(-q)^{-p}(q^2;q^2)_{\infty}(q^{2p+2};q^2)_{\infty}^{-1/2}$, and using the definition of $P_{q^2}^0$ (see Proposition \ref{PropBasEq2}) and the expression (\ref{eqnPfunc}) for the function $P_{q^2}^+$ (see Proposition \ref{PropBasSUq2}), this becomes \begin{eqnarray*} && \underset{v-w=t}{\sum_{v,w\in I_0}} (-q)^{-w}q^{(p-w)(v-w)}\qortPsi{0}{q^{2v-2w+2}}{q^2}{q^{2p-2w+2}} e_{v}\otimes e_{r+p-w}\otimes e_{w}\otimes e_{s-p+v} \\ &&= \underset{v-w=t}{\sum_{v,w\in I_+}}\sum_{k=0}^{\infty} (-q)^{k}q^{vw+(p+k)(v+w)}\frac{(q^{2v+2},q^{2w+2},q^{2k+2};q^2)_{\infty}}{(q^2;q^2)_{\infty}}\qortb{q^{-2w}}{q^{-2v}}{q^{-2p}}{0}{0}{q^2}{q^2} \\ && \qquad \qquad \qquad e_{v-k}\otimes e_{r+p-w+k}\otimes e_{w-k}\otimes e_{s-p+v-k}.\end{eqnarray*} But the coefficients on the right hand side make sense for all $v,w\in I_0$, and are moreover $=0$ when $v$ or $w$ is not in $I_+$. We may hence take the summation on the right over $I_0$, and then a comparison of coefficients shows that we must prove the following summation formula: for all $v,w\in \mathbb{Z}$ and $p\in \mathbb{N}$, \begin{eqnarray*}&&\!\!\!\!\!\!\!\!\!\!(-q)^{-w}q^{(p-w)(v-w)}\qortPsi{0}{q^{2v-2w+2}}{q^2}{q^{2p-2w+2}} = \\ &&\qquad \qquad \qquad\qquad\qquad \sum_{k=-\infty}^{\infty} (-q)^k q^{(v+k)(w+k)+(p+k)(v+w+k)}\frac{(q^{2v+2k+2},q^{2w+2k+2},q^{2k+2};q^2)_{\infty}}{(q^2;q^2)_{\infty}}\\ && \qquad \qquad \qquad \qquad \qquad \qquad \cdot\qortb{q^{-2w-2k}}{q^{-2v-2k}}{q^{-2p}}{0}{0}{q^2}{q^2}.\end{eqnarray*} But, with $t=v-w$, the right hand side can be rewritten as \[\sum_{k=0}^{\infty} (-q)^{k-w}q^{3k^2+2(t+p-w)k+(p-w)t} \frac{(q^{2k-2w+2},q^{2t+2k+2};q^2)_{\infty}}{(q^2;q^2)_k}\qortb{q^{-2k}}{q^{-2t-2k}}{q^{-2p}}{0}{0}{q^2}{q^2} ,\] and the identity then follows by applying Proposition \ref{PropSum3} with respect to $x=q^{-2w}$ and $y=q^{2t}$. This finishes the proof.

\end{proof}

\noindent Combining this result with the work of the previous section, it is now an easy task to finish the proof of Theorem \ref{TheoCoa}.

\begin{proof}[Proof (of Theorem \ref{TheoCoa})] Using the notation introduced before the Theorem, we have shown in the previous sections that all maps $\Delta_{\mu\nu}$ are coassociative, except for those with $\mu\neq \nu$ inside $\{-,0\}$. But take $x\in \mathscr{L}(0,+)$ and $y\in \mathscr{L}(+,-)$. Then from the definition of the $\Delta_{\mu\nu}$, it follows immediately that $\Delta_{0-}(xy) = \Delta_{0+}(x)\Delta_{+-}(y)$, and hence \begin{eqnarray*} (\Delta_{0-}\otimes \iota)\Delta_{0-}(xy) &=& (\Delta_{0+}\otimes \iota)\Delta_{0+}(x)(\Delta_{-0}\otimes \iota)\Delta_{-0}(y) \\ &=& (\iota\otimes \Delta_{0+})\Delta_{0+}(x)(\iota\otimes \Delta_{+-})\Delta_{+-}(y)\\ &=& (\iota\otimes \Delta_{0-})\Delta_{0-}(xy),\end{eqnarray*} where the second identity follows from the results of the previous sections. As the set \[\mathscr{L}(0,+)\cdot \mathscr{L}(+,-)\subseteq \mathscr{L}(0,-)\] is $\sigma$-weakly dense in $\mathscr{L}(0,-)$, we have proven that $\Delta_{0-}$ is coassociative. The coassociativity of $\Delta_{-0}$ then follows immediately by applying the $^*$-operation. This concludes the proof.

\end{proof}

\noindent \emph{Remark:} One can also prove an analogue of Theorem \ref{TheoCoa} if we replace $SU_q(2)$, $\widetilde{E}_q(2)$ and $\widetilde{SU}_q(1,1)$ by the respective $\mathbb{Z}_2$-quotients $SO_q(3)$, $E_q(2)$ and $\widetilde{SU}_q(1,1)/\mathbb{Z}_2\cong O_q(1,2)$ (although this final quantum group may require a different interpretation, see the remark after Definition \ref{DefSUq11}). The reason is that the actions of $SU_q(2)$ on the Podle\'s spheres descend to $SO_q(3)$, so that one may equally well apply our theory of projective corepresentations to $SO_q(3)$ with respect to the standard Podle\'s sphere and the quantum projective plane, as to obtain in this way a 3$\times$3-linking weak von Neumann bialgebra between the mentioned quantum groups. In fact, the only reason why we have worked with their double covers is that these have received more attention in the literature, so that it was easier for us to refer to the known results concerning these quantum groups.\\

\appendix

\section{On the tensor product representations of $\mathscr{L}^{\infty}(SU_q(2))$, $\mathscr{L}^{\infty}(\widetilde{E}_q(2))$ and $\mathscr{L}^{\infty}(\widetilde{SU}_q(1,1))$}\label{ApA}

\noindent Let us comment on the proof of Proposition \ref{PropBasSUq2}, Proposition \ref{DefEq2}, and the equivalence of the Definition \ref{DefSUq11} with the one in \cite{Koe3}. We will use the same notations as in the introduction.\\

\noindent For the case of $SU_q(2)$, there are several references to give. In \cite{Koo1}, the spectral decomposition of the operator $\Delta_+(b_+^*b_+)$ was found. This was given in terms of a basis $\widetilde{F}_{r,s,p,t}^{+}$ of $l^2(\mathbb{Z})\otimes l^2(\mathbb{Z})\otimes \mathscr{H}_+$. (The $r$-$s$-variables were not present in that description, as one worked with a different, irreducible representation of the polynomial function algebra of $SU_q(2)$. But the insertion of the extra $r$-$s$-variables is entirely straightforward, and simply a question of doing some small bookkeeping on shifts.) However, for the action on the $p$-$t$-part to be \emph{exactly} the same as the original representation of $\mathscr{L}^{\infty}(SU_q(2))$ on $\mathscr{H}_+$, one should add some minus-signs: the correct basis then becomes $F_{r,s,p,t}^+ = (-1)^p \widetilde{F}_{r,s,p,t}^+$. That the associated unitary $F_{r,s,p,t}\rightarrow e_{r}\otimes e_s\otimes e_{p}\otimes e_t$ then implements the comultiplication on $\mathscr{L}^{\infty}(SU_q(2))$ was presented in detail in \cite{Gro1} (we would like to thank W. Groenevelt for providing this manuscript). Now the concrete formula for the $F_{r,s,p,t}^+$ is \[F_{r,s,p,t}^+ = \sum_{w\in \mathbb{N}} (-1)^pP_w(q^{2p};q^{2t}\mid q^2) e_{w+t}\otimes e_{r+p-w} \otimes e_w\otimes e_{s-p+t+w},\] where $P_w$ is the \emph{normalized Wall polynomial}: \[P_w(q^{2p};q^{2t}\mid q^2) = (-1)^w q^{(p-w)(t+1)} \frac{(q^{2t+2},q^{2p+2};q^2)_{\infty}^{1/2}(q^{2t+2};q^2)_{w}}{(q^2;q^2)_\infty^{1/2}(q^2;q^2)_{w}^{1/2}} p_{w}(q^{2p};q^{2t},0\mid q^2),\] with $p_{w}(q^{2p};q^{2t},0\mid q^2)= \qorta{q^{-2w}}{0}{q^{2t+2}}{q^2}{q^{2p+2}}$ the ordinary Wall polynomial (of degree $w$ in the variable $q^{2p}$ with parameter $q^{2t}$). One then uses a limit version of the identity III.(3) of \cite{Gas1} to arrive at the expression we used in Proposition \ref{PropBasSUq2} (this was also observed in the remarks following Proposition 3.3 of \cite{Koe3}). On the other hand, using the transformation which appears at the end of section 2 of \cite{Koo1}, we get the expression (\ref{eqnPfunc}) for the functions $P_{q^2}^+$ appearing in Proposition \ref{PropBasSUq2}.\\

\noindent We remark that an implementing unitary for $\Delta_+$ (in fact, a concrete description of the \emph{regular left corerepresentation of $SU_q(2)$}) was also obtained in \cite{Lan1}.\\

\noindent Proposition \ref{PropBasEq2} was proven in the course of Proposition \ref{PropIdLink}. However, we should remark that in \cite{Koe1}, a spectral decomposition of the operator $\Delta_0(b_0^*b_0)$ was obtained. This gives an eigenvector-basis that, as with $\cite{Koo1}$ in the previous paragraph, is almost the same basis as the one we obtain, but modulo the appearance of some extra minus sign. We should also note that \[P_{q^2}^{0}(p;v;w) = (-q)^{p-w}J_{v-w}(q^{p-w};q^2),\] where $J_{\alpha}(z;q^2)$ is the $\,\!_1\varphi_1$-$q$-Bessel function \[J_{\alpha}(z;q^2)=z^{\alpha} \frac{1}{(q^2;q^2)_{\infty}} \qortPsi{0}{q^{2\alpha+2}}{q^2}{q^{2}z^{2}}.\]
\vspace{0.2cm}

\noindent Finally, we must comment on the definition of $\widetilde{SU}_q(1,1)$ we presented. The equivalence with the one in \cite{Koe3} is again straightforward, though there is now an extra element which should be observed.\\

\noindent We first introduce a certain operator $\Omega\in\mathscr{L}^{\infty}(\widetilde{SU}_q(1,1))\bar{\otimes}\mathscr{L}^{\infty}(\widetilde{SU}_q(1,1))$. As $e$ is a self-adjoint unitary inside $\mathscr{L}^{\infty}(\widetilde{SU}_q(1,1))$, we may identify $W^*(e)$ with $C(\mathbb{Z}_2)$, sending $e$ to the identity function. As $e$ is group-like, this will moreover intertwine $\Delta_-$ with the natural comultiplication on $C(\mathbb{Z}_2)$ coming from the group structure on $\mathbb{Z}_2$. Denote then by $\omega \in C(\mathbb{Z}_2)\otimes C(\mathbb{Z}_2)$ the 2-cocycle function \[\omega:\mathbb{Z}_2\times \mathbb{Z}_2\rightarrow S^1: \left\{\begin{array}{lll} \omega(\mu,\nu)=1& \textrm{if} & \mu \textrm{ and } \nu \textrm{ not both }- \\ \omega(\mu,\nu)=-1 &\textrm{if}&\mu=\nu=-,\end{array}\right.\] and denote by $\Omega$ its image in $\mathscr{L}^{\infty}(\widetilde{SU}_q(1,1))\bar{\otimes}\mathscr{L}^{\infty}(\widetilde{SU}_q(1,1))$, using the above identification. More explicitly, we have that $\Omega$ is given on $\mathscr{H}_-\otimes \mathscr{H}_-$ as the operator \[\Omega\,e_{\mathbf{v}}\otimes e_{r}\otimes e_{\mathbf{w}}\otimes e_s = \omega(c(\mathbf{v}),c(\mathbf{w}))e_{\mathbf{v}}\otimes e_{r}\otimes e_{\mathbf{w}}\otimes e_s.\] The 2-cocycle property of $\omega$ will give the following identity for $\Omega$: \[(\Omega_-\otimes 1)(\Delta_-\otimes \iota)(\Omega_-) = (1\otimes \Omega_-)(\iota\otimes \Delta_-)(\Omega_-),\] i.e. $\Omega$ is a unitary $2$-cocycle (\cite{Eno2}) for $\Delta_-$. This implies that $(\mathscr{L}^{\infty}(\widetilde{SU}_q(1,1),\Omega\Delta_-(\,\cdot\,)\Omega^*)$ is again a well-defined von Neumann bi-algebra.\\

\noindent Let us now turn to the definition of $\widetilde{SU}_q(1,1)$ presented in \cite{Koe3}. We first remark that in \cite{Koe3}, one identifies our set $I_-$ with a subset of $\mathbb{R}$ by the correspondence $p_{\pm} \leftrightarrow \pm q^{-p}$. Denote then by $\mathcal{F}_{r,s,m,\pm q^{-p}}$ the basis elements introduced in Definition 3.6 of \cite{Koe3}. One computes, using Result 6.4 and Proposition 6.6 of \cite{Koe3}, that the following relation holds between the vector bases $\mathcal{F}$ and $\xi^-$: \[\mathcal{F}_{r,s,m,\pm q^{-p}} = \Omega \xi_{r,s,\mathbf{p},m}^-.\] This implies that, if we denote by $\widetilde{\Delta}_-$ the comultiplication as introduced in \cite{Koe3} by means of the vector basis $\mathcal{F}$, we have the following relation: \[ \widetilde{\Delta}_-(x) = \Omega\Delta_-(x)\Omega^*,\qquad \textrm{for all }x\in \mathscr{L}^{\infty}(\widetilde{SU}_q(1,1)).\]

\noindent But in fact, $(\mathscr{L}^{\infty}(\widetilde{SU}_q(1,1)),\Omega\Delta_-(\,\cdot\,)\Omega^*)$ is just an isomorphic copy of $(\mathscr{L}^{\infty}(\widetilde{SU}_q(1,1),\Delta_-(\,\cdot\,))$. To see this, we remark that, as we have taken $\omega$ as an \emph{$S^1$-valued} 2-cocycle on $\mathbb{Z}_2$, it is a coboundary: we have $\omega(\mu,\nu) = \frac{f(\mu\nu)}{f(\mu)f(\nu)}$ with $f(+)=1$ and $f(-)=i$. Denote then \[u : \mathscr{H}_-\rightarrow \mathscr{H}_-: e_{\mathbf{p}}\otimes e_r\rightarrow f(c(\mathbf{p}))\,e_{\mathbf{p}}\otimes e_r,\] and denote \[\phi: \mathscr{L}^{\infty}(\widetilde{SU}_q(1,1))\rightarrow \mathscr{L}^{\infty}(\widetilde{SU}_q(1,1)): x\rightarrow uxu^*.\] Then one gets that $\Omega = (u^*\otimes u^*)\Delta_-(u)$, and so \[\Delta_-(\phi(x)) = (\phi\otimes \phi)(\Omega\Delta_-(x)\Omega^*), \qquad \textrm{for all }x\in \mathscr{L}^{\infty}(\widetilde{SU}_q(1,1)).\] This thus proves that $(\mathscr{L}^{\infty}(\widetilde{SU}_q(1,1)),\Delta_-)$ and $(\mathscr{L}^{\infty}(\widetilde{SU}_q(1,1)),\widetilde{\Delta}_-)$ are isomorphic von Neumann bialgebras, and justifies our use of the basis $\xi^-$.

\section{Spectral decomposition of $\widetilde{\alpha}_+(W)$}\label{ApB}

\noindent In this appendix, we will determine the spectral decomposition of the action of $SU_q(2)$ on the equatorial Podle\'{s} sphere and on the quantum projective plane. We are not aware of this concrete computation having been carried out explicitly in the literature, but the method is quite standard (and could probably be shortened somewhat).\\

\noindent We will use the notation as introduced in section \ref{SecPod}.

\begin{Prop} The eigenspace of the eigenvalue $\pm 1$ of $\widetilde{\alpha}_+(W)$ has a basis of orthonormal vectors $\xi_{0,\pm}^{(t,p)}$, where $t\in \mathbb{Z}$ and $p\in J^{(\pm)}$, these vectors being given by the formula \[\xi_{0,\pm}^{(t,p)} = \sum_{n=0}^{\infty} (\pm 1)^n  q^{\frac{n(n-1)}{2}}\frac{(\mp q^{2p+2};q^2)_{n}^{1/2}(\pm q^{2p+2n+2};q^2)_{\infty}^{1/2}}{(-1;q^2)_{\infty}^{1/2}(q^2;q^2)_n^{1/2}} \, e_{n}\otimes e_{t-n} \otimes e_{p+n}.\]
\end{Prop}

\begin{proof}

\noindent Denote, for $t,p\in \mathbb{Z}$, \[\mathscr{K}^{t,p} = \langle \{ e_n\otimes e_m \otimes e_k \mid m+n = t, k-n=p\} \rangle  \subseteq l^2(\mathbb{N})\otimes l^2(\mathbb{Z})\otimes l^2(\mathbb{N}),\] where $\langle S \rangle$ denotes the linear span of a set $S$. Then, since for $n,k\in \mathbb{N}$ and $m\in \mathbb{Z}$, we have, by the definition of $\widetilde{\alpha}_+(W)$ in Definition \ref{DefCoacProj}, \begin{eqnarray*} \widetilde{\alpha}_+(W)(e_n\otimes e_m\otimes e_k) &=& q^n(1-q^{2n+2})^{1/2} (1-q^{4k+4})^{1/2} e_{n+1} \otimes e_{m-1}\otimes e_{k+1} \\ && \quad + (1-(1+q^2)q^{2n})q^{2k} e_n\otimes e_m\otimes e_k \\ && \quad + q^{n-1} (1-q^{2n})^{1/2} (1-q^{4k})^{1/2} e_{n-1} \otimes e_{m+1} \otimes e_{k-1},\end{eqnarray*} we see that $\widetilde{\alpha}_+(W)$ restricts to each $\mathscr{K}^{t,p}$, and determines there a Jacobi matrix. \\

\noindent We make a distinction in the analysis between the case $p\geq 0$ and $p<0$.\\

\noindent First we consider the case $p\geq 0$. Then we have a unitary map $l^2(\mathbb{N})\rightarrow \mathscr{K}^{t,p}$ such that $e_n\rightarrow e_n\otimes e_{t-n} \otimes e_{p+n}$. Under this identification, the restriction of $\widetilde{\alpha}_+(W)$ becomes the Jacobi matrix  $W^{(p)}$ with \begin{eqnarray*} W^{(p)} e_n &=& q^n(1-q^{2n+2})^{1/2} (1-q^{4(p+n)+4})^{1/2} e_{n+1} \\ && \quad + (1-(1+q^2)q^{2n})q^{2(p+n)} e_n \\ && \quad + q^{n-1} (1-q^{2n})^{1/2} (1-q^{4(p+n)})^{1/2} e_{n-1}.\end{eqnarray*} Each eigenvalue then arises with multiplicity one, and an eigenvector at eigenvalue $x$ is given by $\eta_{x}^{(p)} = \sum f_n^{(p)}(x)e_n$, where $f_n^{(p)}$ is the sequence of functions satisfying $f^{(p)}_{-1} = 0$, $f^{(p)}_0 = 1$ and \[q^{n}(1-q^{2n+2})^{1/2} (1-q^{4p+4n+4})^{1/2} f_{n+1}^{(p)}(x) + (1-(1+q^2)q^{2n})q^{2p+2n} f_n^{(p)}(x) \]\[\qquad \qquad +q^{n-1} (1-q^{2n})^{1/2} (1-q^{4p+4n})^{1/2} f_{n-1}^{(p)}(x) = xf_n^{(p)}(x).\]

\noindent Denote \[f_n^{(p)}(x) = (-1)^n q^{-\frac{n^2}{2} - \frac{4p+3}{2}n} \frac{(q^4;q^4)_{p+n}^{1/2}}{(q^4;q^4)_{p}^{1/2} (q^{2};q^2)_{n}^{1/2}} P_n^{(p)}(-q^{2p+2}x),\] for some other function $P_n^{(p)}$. Then $P_n^{(p)}$ should itself satisfy the recurrence relation \[\!\!\!\!\!\!\!\!\!\!\!\!\!\!\!\!\!\!\!\!(1-q^{4(p+n+1)})P_{n+1}^{(p)}(x) - ((1-(1+q^2)q^{2n})q^{2n+4p+2}+1)P_{n}^{(p)}(x)\]\[\qquad \qquad \qquad \qquad \qquad \qquad + q^{2n+4p+2}(1-q^{2n})P_{n-1}^{(p)}(x) = (x-1)P_{n}^{(p)}(x).\] Hence, for example by \cite{Koek1}, section 3.11, we see that $P_{n}^{(p)}(x) = P_n(x;q^{2p},-q^{2p};q^2)$, where $P_n(x;a,b;q)$ denotes the big $q$-Laguerre polynomial \[P_n(x;a,b;q) = \qortb{q^{-n}}{0}{x}{qa}{qb}{q}{q}.\] So \[f_n^{(p)}(x) = (-1)^n q^{-\frac{n^2}{2} - \frac{4p+3}{2}n} \frac{(q^4;q^4)_{p+n}^{1/2}}{(q^4;q^4)_{p}^{1/2} (q^{2};q^2)_{n}^{1/2}}\,\qortb{q^{-2n}}{0}{-q^{2p+2}x}{q^{2p+2}}{-q^{2p+2}}{q^2}{q^2}.\]

\noindent Now we know that the eigenspace $\mathscr{H}_{\pm 1}$ for the eigenvalue $\pm 1$ of $\widetilde{\alpha}_+(W)$ will be spanned by the $\mathscr{H}_{\pm 1} \cap \mathscr{K}^{t,p}$. We now show that each $W^{(p)}$ has eigenvalues $1$ and $-1$.\\

\noindent We first remark that \[ \qortb{q^{-2n}}{0}{\mp q^{2p+2}}{q^{2p+2}}{-q^{2p+2}}{q^2}{q^2} = \qorta{q^{-2n}}{0}{\pm q^{2p+2}}{q^2}{q^2}.\] Then, using a limit form of the $q$-Vandermonde formula (\cite{Gas1}, Equation 1.5.3), we get \[\qorta{q^{-2n}}{0}{\pm q^{2p+2}}{q^2}{q^2}=(\mp 1)^n q^{n(n-1)} \frac{q^{2(p+1)n}}{(\pm q^{2p+2};q^2)_n}.\] Hence, using the formula $(a^2;q^2)_n = (a;q)_n(-a;q)_n$, we get \begin{eqnarray*} f_{n}^{(p)}(\pm 1) &=& (\pm 1)^n q^{\frac{n(n-1)}{2}} \frac{(\mp q^2;q^2)_{p+n}^{1/2}}{(q^2;q^2)_n^{1/2} (\mp q^2;q^2)_p^{1/2} (\pm q^{2p+2};q^2)_n^{1/2}}\\ &=& (\pm 1)^n q^{\frac{n(n-1)}{2}}\frac{(\mp q^{2p+2};q^2)_{n}^{1/2}}{(q^2;q^2)_n^{1/2}(\pm q^{2p+2};q^2)_n^{1/2}}.\end{eqnarray*} So, by using a limit form of Heine's summation formula (\cite{Gas1}, Equation II.5), \begin{eqnarray*} \sum_{n=0}^{\infty} |f_n^{(p)}(\pm 1)|^2 &=& \qortc{\mp q^{2p+2}}{\pm q^{2p+2}}{q^2}{-1} \\ &=& \frac{(-1;q^2)_{\infty}}{(\pm q^{2p+2};q^2)_{\infty}}.\end{eqnarray*}

\noindent Hence $\pm 1$ appears as an eigenvalue of $W^{(p)}$ with eigenvector $\eta_{\pm 1}^{(p)}\in \mathscr{K}^{t,p}$. We then find that the intersection of the eigenspace $\mathscr{H}_{\pm 1}$ with the closed linear span of the $\mathscr{K}^{t,p}$ with $p\geq0$ has an orthonormal basis consisting of the vectors \[ \xi_{0,\pm}^{(t,p)} = \sum_{n=0}^{\infty} (\pm 1)^n q^{\frac{n(n-1)}{2}}\frac{(\mp q^{2p+2};q^2)_{n}^{1/2}(\pm q^{2p+2n+2};q^2)_{\infty}^{1/2}}{(-1;q^2)_{\infty}^{1/2}(q^2;q^2)_n^{1/2}} \, e_{n}\otimes e_{t-n} \otimes e_{p+n}.\]

\vspace{0.4cm}

\noindent We move on to the case $p<0$. Now we have a unitary map $l^2(\mathbb{N})\rightarrow \mathscr{K}^{t,p}$ such that $e_k$ corresponds to $e_{k-p} \otimes e_{t-k+p} \otimes e_{k}$. Under this identification, the restriction of $\widetilde{\alpha}_+(W)$ becomes the Jacobi matrix  $W^{(p)}$ with \begin{eqnarray*} W^{(p)} e_k &=& q^{k-p}(1-q^{2k-2p+2})^{1/2} (1-q^{4k+4})^{1/2} e_{k+1} \\ && \quad + (1-(1+q^2)q^{2(k-p)})q^{2k} e_k \\ && \quad + q^{k-p-1} (1-q^{2k-2p})^{1/2} (1-q^{4k})^{1/2} e_{k-1}.\end{eqnarray*} Again, each eigenvalue arises with multiplicity one, and an eigenvector at eigenvalue $x$ is given by $\eta_{x}^{(p)} = \sum f_k^{(p)}(x)e_k$, where now $f_k^{(p)}$ is the sequence of functions satisfying $f^{(p)}_{-1} = 0$, $f^{(p)}_0 = 1$ and \[ q^{k-p} (1-q^{2k-2p+2})^{1/2} (1-q^{4k+4})^{1/2} f_{k+1}^{(p)}(x) + (1-(1+q^2)q^{2(k-p)})q^{2k} f_{k}^{(p)}(x) \]\[ \quad + q^{k-p-1}(1-q^{2k-2p})^{1/2} (1-q^{4k})^{1/2} f_{k-1}^{(p)}(x) = xf_k^{(p)}(x).\] Now put \[f_{k}^{(p)}(x) = (-1)^k q^{-\frac{k^2}{2} - \frac{(-2p+3)k}{2}}\, \frac{(q^{-2p+2};q^2)_k^{1/2}(-q^2;q^2)_k^{1/2}}{(q^2;q^2)_k^{1/2}} P_k^{(p)}(-q^2x),\] for some function $P_k^{(p)}$. This $P_k^{(p)}$ should then satisfy the recursion formula \[(1-q^{2k-2p+2})(1+q^{2k+2}) P_{k+1}^{(p)}(x) - ((1-(1+q^2)q^{2k-2p})q^{2k+2}+1)P_{k}^{(p)}(x)\]\[\qquad + q^{2k-2p+2}(1-q^{2k})P_{k-1}^{(p)}(x) = (x-1)P_{k}^{(p)}(x).\] So $P_k^{(p)}$ is again a big $q$-Laguerre polynomial, namely $P_k^{(p)}(x) = P_k(x;q^{-2p},-1;q^2)$. We then have \[f_{k}^{(p)}(x) = (-1)^k q^{-\frac{k^2}{2} - \frac{(-2p+3)k}{2}}\, \frac{(q^{-2p+2};q^2)_k^{1/2}(-q^2;q^2)_k^{1/2}}{(q^2;q^2)_k^{1/2}} \qortb{q^{-2k}}{0}{-q^2x}{q^{-2p+2}}{-q^2}{q^2}{q^2}.\] We now show that $W^{(p)}$ has $1$ in its spectrum, but not $-1$.\\

\noindent We first show that $1$ is in the spectrum. Using again the $q$-Vandermonde formula, we have \begin{eqnarray*} \qortb{q^{-2k}}{0}{-q^2}{q^{-2p+2}}{-q^2}{q^2}{q^2} &=& \qorta{q^{-2k}}{0}{q^{-2p+2}}{q^2}{q^2} \\ &=& (-1)^k q^{k(k-1)} \frac{q^{-2pk+2k}}{(q^{-2p+2};q^2)_k}.\end{eqnarray*} Then \[f_k^{(p)}(1) = q^{\frac{k(k-1)}{2}} q^{-pk} \frac{(-q^2;q^2)_k^{1/2}}{(q^2;q^2)_k^{1/2}(q^{-2p+2};q^2)_k^{1/2}}.\] So, using again the limit form of Heines summation formula, we get \begin{eqnarray*} \sum_{k=0}^{\infty} |f_k^{(p)}(1)|^2 &=& \qortc{-q^2}{q^{-2p+2}}{q^2}{-q^{-2p}}\\ &=& \frac{(-q^{-2p};q^2)_{\infty}}{(q^{-2p+2};q^2)_{\infty}}.\end{eqnarray*} Hence the formal sum $\eta_{1}^{(p)}$ gives in fact a well-defined element inside $\mathscr{K}^{t,p}$, and we find that the intersection of the eigenspace $\mathscr{H}_{1}$ with the closed linear span of the $\mathscr{K}^{t,p}$ with $p<0$ has an orthonormal basis consisting of the vectors \[ \xi_{0,+}^{(t,p)} = \sum_{k=0}^{\infty} q^{\frac{k(k-1)}{2}} q^{-pk} \frac{(-q^2;q^2)_k^{1/2}(q^{2k-2p+2};q^2)_{\infty}^{1/2}}{(-q^{-2p};q^2)_{\infty}^{1/2}(q^2;q^2)_k^{1/2}} \,e_{k-p}\otimes e_{t-k+p} \otimes e_k.\] Now for $n\geq -p$, we have \[ \frac{(-q^2;q^2)_{n+p}^{1/2} (q^{2n+2};q^2)_{\infty}^{1/2}}{(-q^{-2p};q^2)_{\infty}^{1/2}(q^2;q^2)_{n+p}^{1/2}} = \frac{(-q^2;q^2)_{\infty}^{1/2}}{(-q^{2n+2p+2};q^2)_{\infty}^{1/2} (-q^{-2p};q^2)_{\infty}^{1/2}}\cdot \frac{(q^{2n+2p+2};q^2)_{\infty}^{1/2}}{(q^2;q^2)_n^{1/2}},\] and \[ \frac{(-q^2;q^2)_{\infty}^{1/2}}{(-q^{2n+2p+2};q^2)_{\infty}^{1/2} (-q^{-2p};q^2)_{\infty}^{1/2}} = q^{\frac{(-p)(-p-1)}{2}} \frac{(-q^{2p+2};q^2)_n^{1/2}}{(-1;q^2)_{\infty}^{1/2}}.\] Hence we can also write \[ \xi_{0,+}^{(t,p)} = \sum_{n=0}^{\infty}  q^{\frac{n(n-1)}{2}}\frac{(- q^{2p+2};q^2)_{n}^{1/2}( q^{2p+2n+2};q^2)_{\infty}^{1/2}}{(-1;q^2)_{\infty}^{1/2}(q^2;q^2)_n^{1/2}} \, e_{n}\otimes e_{t-n} \otimes e_{p+n},\] which is precisely the same formula as for the $p\geq 0$ case. Indeed, we could also have checked directly that this is also an eigenvector of norm 1 for the eigenvalue $1$.\\

\noindent Now we show that $W^{(p)}$ has no eigenvector for $-1$ when $p<0$. First remark that, using the transformation formulas (III.5) and (III.1) of \cite{Gas1}, we can rewrite, for $-p>0$, \begin{eqnarray*}
&&\!\!\!\!\!\!\!\!\!\!\!\!\qortb{q^{-2k}}{0}{q^2}{\!\!q^{-2p+2}\,\,}{\,\,-q^2}{q^2}{q^2} \\&& = (-q^{-2k};q^2)_k^{-1}\, \qorta{q^{-2k}}{q^{-2p}}{q^{-2p+2}}{q^2}{-q^2}\\ &&= \frac{(q^{-2p};q^2)_{\infty}(-q^{-2k+2};q^2)_{\infty}}{(-q^{-2k};q^{2})_k(q^{-2p+2};q^2)_{\infty}(-q^2;q^2)_{\infty}}
\,\qorta{q^2}{-q^2}{-q^{-2k+2}}{q^2}{q^{-2p}}.\end{eqnarray*} But clearly \[ 1 \leq \qorta{q^2}{-q^2}{-q^{-2k+2}}{q^2}{q^{-2p}} \leq \qorta{q^2}{-q^2}{-q^2}{q^2}{q^{-2p}}\] for all $-p> 0$. So to see if $\sum_{k=0}^{\infty} |f_k^{(p)}(-1)|^2 =+\infty$, we have to check if \[\sum_{k=0}^{\infty} q^{-k^2-(-2p+3)k} \, \frac{(q^{-2p+2};q^2)_k(-q^2;q^2)_k(-q^{-2k+2};q^2)_{\infty}^{2}}{(q^2;q^2)_k (-q^{-2k};q^2)_k^{2}} = +\infty.\] Leaving out more non-essential factors regarding the convergence, the sum simplifies to \[\sum_{k=0}^{\infty}   q^{-k^2-(-2p+3)k} (1+q^{-2k})^{-2},\] which is clearly divergent.\\

\end{proof}

\noindent It is now not so hard to find the complete spectrum of $\widetilde{\alpha}_+(W)$ by using the operators $\widetilde{\alpha}_+(Y)$ and its adjoint, as given in Definition \ref{DefCoacProj}.

\begin{proof}[Proof (of Proposition \ref{PropEigvect})]
\noindent The fact that $\widetilde{\alpha}_+(W)$ will have its spectrum inside $\{\pm q^{2n},0\mid n\in \mathbb{N}\}$, with 0 not in the point-spectrum, is immediate.\\

\noindent Using the same notation as in the previous proposition, denote, for $r\in \mathbb{N}$ and $t,p\in \mathbb{Z}$, by $\mathscr{K}_{r,\pm}^{t,p}$ the eigenspace of the eigenvector $\pm q^{2r}$ inside $\mathscr{K}^{t,p}$ (which could be zero-dimensional). Then for $r>0$, we have, by using the commutation relations between $W,Y$ and $Y^*$, \[\widetilde{\alpha}_+(Y) \mathscr{K}_{r,\pm}^{t,p} = \mathscr{K}_{r-1,\pm}^{t-2,p+1},\]\[ \widetilde{\alpha}_+(Y)^* \mathscr{K}_{r,\pm}^{t,p} = \mathscr{K}_{r+1,\pm}^{t+2,p-1}.\]

\noindent Denote, for $r\in \mathbb{N}$ and $t,p\in \mathbb{Z}$, by $\eta_{r,\pm}^{(t,p)}$ the formal eigenvectors which we denoted as $\eta_{\pm q^{2r}}^{(p)} \in \mathscr{K}^{t,p}$ in the previous Proposition. Then by combining the above remark concerning $\widetilde{\alpha}_+(Y)$ and its adjoint with the previous Proposition, we see that in the $+$-case, all values of $p$ and $r$ correspond to actual eigenvectors, while in the $-$-case, we have the restriction $r+p\geq 0$. Moreover, the linear span of the $\eta_{r,\pm}^{(t,p)}$, where $r,\pm,t$ and $p$ range over all admissible values, will then be a dense subspace of $l^2(\mathbb{N})\otimes l^2(\mathbb{Z})\otimes l^2(\mathbb{N})$.\\

\noindent Now we have, by comparing leading coefficients in the standard basis decomposition, for $p>0$ that \begin{equation}\label{EqActY} \widetilde{\alpha}_+(Y^*) \eta_{r,\pm}^{(t,p)} = -q(1-q^{4p})^{1/2} \, \eta_{r+1,\pm}^{(t+2,p-1)},\end{equation} while for $p\leq0$, we have  \begin{equation}\label{EqActY2} \widetilde{\alpha}_+(Y^*) \eta_{r,\pm}^{(t,p)} = -\frac{(1\pm q^{2r+2})}{(1-q^{-2p+2})^{1/2}}q^{-p+1} \, \eta_{r+1,\pm}^{(t+2,p-1)}.\end{equation} Since $YY^* = 1-q^4W^2$ by an easy computation, we then have, taking norms of both sides, that, for $p>0$, \[ (1-q^{4r+4})\|\eta_{r,\pm}^{(t,p)}\|^2 = q^2(1-q^{4p})\|\eta_{r+1,\pm}^{(t+2,p-1)}\|^2,\] and for $p\leq0$, \[(1-q^{4r+4}) \| \eta_{r,\pm}^{(t,p)}\|^2 =  \frac{(1\pm q^{2r+2})^2}{(1-q^{-2p+2})}q^{-2p+2} \, \|\eta_{r+1,\pm}^{(t+2,p-1)}\|^2.\]

\noindent By induction, we then find the following formulas for the norm squared of the $\eta$-vectors: if $p\geq 0$, we have \[\|\eta_{r,\pm}^{(t,p)}\|^2 = q^{-2r} \frac{(q^4;q^4)_r(-1;q^2)_{\infty}}{(q^{4p+4};q^4)_r(\pm q^{2p+2r+2};q^2)_{\infty}},\] if $r+p\geq 0>p$, we have \[\|\eta_{r,\pm}^{(t,p)}\|^2 = q^{-p^2-p-2r}  \frac{(q^2;q^2)_{-p} (\mp q^{2r+2};q^2)_{\infty}(-1;q^2)_{\infty}}{(\pm q^{2r+2};q^2)_{\infty}(\mp q^{2r+2p+2};q^2)_{\infty}},\] and for $0>r+p>0$ we have \[ \|\eta_{r,+}^{(t,p)}\|^2 = q^{r^2-r+2rp} \frac{(q^2;q^2)_r(q^2;q^2)_{-p}(-q^{-2r-2p};q^2)_{\infty}}{(-q^2;q^2)_r(q^2;q^2)_{-p-r} (q^{-2r-2p+2};q^2)_{\infty}}.\]

\noindent Denote by $\xi_{r,\pm}^{(t,p)} = \frac{(\mp 1)^r}{\|\eta_{r,\pm}^{(t,p)}\|}\, \eta_{r,\pm}^{(t,p)}$ the normalized eigenvectors for $\widetilde{\alpha}_+(W)$ (the reason for the extra sign factor $(\mp 1)^r$ is that Proposition \ref{PropIdenG} would hold). Then, from the expressions for $\eta_{r,\pm}^{(t,p)}$ in the previous Proposition, together with the above formulas for the norm, we have concrete expressions for the $\xi_{r,\pm}^{(t,p)}$ in terms of $\!\,_3\varphi_2$-functions. We next want to write these in a different form.\\

\noindent Using, for $p\geq0$, the transformation formula \begin{eqnarray*}\qortb{q^{-2n}}{0}{\mp q^{2p+2r+2}}{\pm q^{2p+2}}{\mp q^{2p+2}}{q^2}{q^2} &=& (-1)^n \frac{(\mp q^{2p+2r+2};q^2)_n}{(\mp q^{-2p-2n};q^2)_n(\pm q^{2p+2};q^2)_n}\\&&\qquad \cdot\qortb{q^{-2n}}{q^{-2r}}{\pm q^{-2p-2n}}{\mp q^{-2p-2n-2r}}{0}{q^2}{q^2},\end{eqnarray*} which is a combination of the identities (III.5) and (III.6) of \cite{Gas1}, we find, after some simplifications, that the $\xi_{r,\pm}^{(t,p)}$ for $p\geq 0$ satisfy the formula in the statement of the Proposition.\\

\noindent In case $p<0$, we can use the same identities to get, for $k=n+p\geq 0$, the transformation formula \begin{eqnarray*} \qortb{q^{-2k}}{0}{\mp q^{2r+2}}{q^{-2p+2}}{-q^2}{q^2}{q^2} &=&(\pm 1)^{k-p} (- q^{-2p})^{k}\frac{(\mp q^{2r+2};q^2)_{k}}{(-q^{-2k};q^2)_k (q^{-2p+2};q^2)_{k}}\\&&\qquad \cdot
\qortb{q^{-2k+2p}}{q^{-2r}}{\pm q^{-2k}}{\mp q^{-2k-2r}}{0}{q^2}{q^2},\end{eqnarray*} from which, again after some simplifications, we see that the formula for $\xi_{r,\pm}^{(t,p)}$ in the Proposition is still valid in this case (where we note again that, in the case of negative eigenvalues, the $r,p$ are restricted by the condition $r+p\geq0$).\end{proof}

\section{Some summation formulas for basic hypergeometric functions}\label{ApC}

\begin{Prop}\label{PropSum1} Suppose $x,y\in \mathbb{C}_0$ and $p\in \mathbb{N}$. For $w\in \mathbb{N}$, write \[g_w^{(p)}(x,y) = (-1)^w q^{2w^2+2pw}(\frac{x}{y})^{w} \frac{(q^{2w+2}\frac{x}{y};q^2)_{\infty}(xq^{2w+2};q^2)_{\infty}(-q^{2}x;q^2)_w}{(q^2;q^2)_w}.\] Then \[\sum_{w=0}^{\infty} g_w^{(p)}(x,y) \qortb{q^{-2w}}{q^{-2w}\frac{y}{x}}{q^{-2p}}{0}{0}{q^2}{q^2} = (-1;q^2)_p \qortPsi{-\frac{1}{y}}{q^2\frac{x}{y}}{q^2}{q^{2p+2}x}.\]
\end{Prop}

\begin{proof} As the $\!\,_3\varphi_2$-term is a polynomial of fixed degree $p$, it is easy to see that the double sum on the left hand side is absolutely convergent, so that we can change the order of summation later on. Moreover, both sides are analytic in each variable $x$ and $y$, so we can also restrict our attention to the situation $x,y\notin \pm q^{2\mathbb{Z}}$.\\

\noindent Then the identity can be rewritten in the form \[\sum_{w=0}^{\infty} (-1)^w q^{2w^2+2pw} (\frac{x}{y})^w  \frac{(-q^{2}x;q^2)_w}{(q^{2}x;q^2)_{w}(q^{2}\frac{x}{y};q^2)_{w}(q^2;q^2)_w} \qortb{q^{-2w}}{q^{2w}\frac{y}{x}}{q^{-2p}}{0}{0}{q^2}{q^2} \]\[= \frac{(-1;q^2)_p}{(q^{2}\frac{x}{y};q^2)_{\infty}(q^{2}x;q^2)_{\infty}} \qortPsi{- \frac{1}{y}}{q^{2}\frac{x}{y}}{q^2}{q^{2p+2}x}.\] \\

\noindent Now on the left hand side, we can expand the $\,\!_3\varphi_2$-term as a sum ranging from $0$ to $w$. Changing the order of summation, we can write the entire expression as a double sum over $l:0\rightarrow \infty$ and $w:l\rightarrow \infty$. Then changing the variable $w$ to $w-l$, and reversing the order of summation again, the above left hand side expression can be simplified to \[\sum_{w=0}^{\infty} (-1)^w q^{2w^2+2pw}(\frac{x}{y})^w \frac{(-q^{2}x;q^2)_w}{(q^{2}\frac{x}{y};q^2)_w(q^{2}x;q^2)_w(q^2;q^2)_w}\sum_{l=0}^{\infty} (-1)^l q^{2lp} \frac{(q^{-2p};q^2)_l (-q^{2w+2}x;q^2)_l}{(q^{2w+2}x;q^2)_l (q^2;q^2)_l}.\] Now the sum over $l$ is just $\qorta{q^{-2p}}{-q^{2w+2}x}{q^{2w+2}x}{q^2}{-q^{2p}}$. By the $q$-Vandermonde formula (\cite{Gas1}, (1.5.2)), this equals $\frac{(-1;q^2)_p}{( q^{2w+2}x;q^2)_p}$. So our left hand side expression becomes \begin{equation}\label{eqn3} \frac{(-1;q^2)_p}{(q^{2}x;q^2)_p}\,\sum_{w=0}^{\infty} (-1)^w q^{2w^2+2pw}(\frac{x}{y})^w \frac{(- q^{2}x;q^2)_w}{(q^{2}\frac{x}{y};q^2)_w(xq^{2p+2};q^2)_w(q^2;q^2)_w}.\end{equation} Now the sum over $w$ can be written as \[\lim_{c\rightarrow 0} \qortd{-q^{2}x}{c^{-1}}{ q^{2p+2}x}{q^{2}\frac{x}{y}}{q^2}{-cq^{2p+2}\frac{x}{y}}.\] Using (III.4) of \cite{Gas1}, this becomes \[\frac{(- \frac{1}{y};q^2)_{\infty}}{(q^{2}\frac{x}{y};q^2)_{\infty}}\qorta{- q^{2}x}{0}{q^{2p+2}x}{q^2}{-\frac{1}{y}},\] which by \cite{Gas1}, (III.1) can be transformed into \[\frac{1}{(q^{2p+2}x,q^{2}\frac{x}{y};q^2)_{\infty}}\, \qortPsi{- \frac{1}{y}}{q^{2}\frac{x}{y}}{q^2}{ q^{2p+2}x}.\] Then plugging this back into the (\ref{eqn3}), we find the identity we were after.\\
\end{proof}

\begin{Prop}\label{PropSum2} Suppose $p\in \mathbb{N}$ and $x,y\in \mathbb{C}_0$ with $x\notin -q^{-2\mathbb{N}-2p-2}$. For $w\in \mathbb{N}$, write \[g_w^{(p,\pm)}(x,y)=(\mp 1)^w q^{2w^2}(\frac{x}{y})^w\frac{(q^{2w+2}\frac{x}{y};q^2)_{\infty}(\pm q^{2w+2}x;q^2)_{\infty}(-q^{2p+2}x;q^2)_w}{(q^2;q^2)_w}.\] Then \begin{eqnarray*} \sum_{w=0}^{\infty} g_w^{(p,\pm)}(x,y)\qortb{q^{-2w}}{ q^{-2w}\frac{1}{x}}{q^{-2p}}{- q^{-2p-2w}\frac{1}{x}}{0}{q^2}{q^2} &=&  \qortPsi{\mp \frac{1}{y}}{\pm q^{2p+2}\frac{x}{y}}{q^2}{\pm  q^{2}x}.\end{eqnarray*}

\end{Prop}

\begin{proof}
\noindent As for the previous Proposition, both sides of the identity are analytic in $x$ and $y$, so we may restrict our attention to the case $x,y\notin \pm q^{2\mathbb{Z}}$. We will then only give the proof for the $+$-case, as the $-$-case is completely similar.\\

\noindent We again expand the $\,\!_3\varphi_2$-factor on the left hand side as a sum over the variable $l:0\rightarrow w$, change the order of $w$ and $l$, replace $w$ by the variable $w-l$, and again change the order of summation. Then we obtain that the expression on the left hand side can be simplified to \[ \sum_{w=0} (-1)^w q^{2w^2}(\frac{x}{y})^w \frac{(-q^{2p+2}x;q^2)_w}{(q^{2}\frac{x}{y};q^2)_w(q^{2}x;q^2)_w(q^2;q^2)_w}  \sum_{l=0}^{\infty} (-1)^l q^{l^2+l}q^{2l(w+p)}(\frac{x}{y})^l \frac{(q^{-2p};q^2)_l}{(q^{2w+2}\frac{x}{y};q^2)_l(q^2;q^2)_l}.\]  The sum over $l$ thus equals $\qortc{q^{-2p}}{q^{2w+2}\frac{x}{y}}{q^2}{q^{2w+2p+2}\frac{x}{y}}$, which, by the limit form of Heines summation formula (\cite{Gas1}, (II.5)) can be reduced to $\frac{(q^{2w+2p+2}\frac{x}{y};q^2)_{\infty}}{(q^{2w+2}\frac{x}{y};q^2)_{\infty}}$. Then the sum over $w$ can be rewritten as \[ (q^{2p+2}\frac{x}{y};q^2)_{\infty} \lim_{c\rightarrow 0} \qortd{- q^{2p+2}x}{c^{-1}}{q^{2}x}{q^{2p+2}\frac{x}{y}}{q^2}{-cq^{2}\frac{x}{y}},\] which, by using again \cite{Gas1}, (III.4) and (III.1), reduces to the right hand side of the identity we wanted to prove.

\end{proof}

\begin{Prop}\label{PropSum3} Suppose $x,y\in \mathbb{C}$ and $p\in \mathbb{N}$. For $k\in \mathbb{N}$, write \[f_k^{(p)} = (-q)^{k}q^{3k^2+2kp}(xy)^{k} \frac{(q^{2k+2}x,q^{2k+2}y;q^2)_{\infty}}{(q^2;q^2)_k}.\] Then  \[ \sum_{k=0}^{\infty} f_k^{(p)}(x,y)\qortb{q^{-2k}}{q^{-2k}\frac{1}{y}}{q^{-2p}}{0}{0}{q^2}{q^2}=   \qortPsi{0}{q^{2}y}{q^2}{q^{2p+2}x}.\]

\end{Prop}

\begin{proof} Both sides of the identity again depend analytically on the variables $x$ and $y$, so we may suppose that $x,y\notin \{0\}\cup q^{2\mathbb{Z}}$. \\

\noindent As in the previous Propositions, we again expand the $\,\!_3\varphi_2$-factor on the left hand side as a sum over the variable $l:0\rightarrow k$, change the order of $k$ and $l$, replace $k$ by the variable $k-l$, and again change the order of summation. Then we obtain that the expression on the left hand side can be simplified to \[\sum_{k=0}^{\infty} (-q)^{k} q^{3k^2+2pk}(xy)^k \frac{(q^{2k+2}x,q^{2k+2}y;q^2)_{\infty}}{(q^2;q^2)_k} \qortc{q^{-2p}}{q^{2k+2}x}{q^2}{q^{2k+2p+2}x}.\] But the $\,\!_1\varphi_1$-expression can be simplified to $\frac{(q^{2k+2p+2}x;q^2)_{\infty}}{(q^{2k+2}x;q^2)_{\infty}}$ by the limit version of Heines summation formula (equation (1.5.1) in \cite{Gas1}). The remaining sum over $k$ can then be shown to equal precisely $\qortPsi{0}{q^{2}y}{q^2}{q^{2p+2}x}$ by a (double) limit version of Jackson's transformation formula (equation (1.5.4) in \cite{Gas1}). This concludes the proof.

\end{proof}

\vspace{0.2cm}

\noindent \emph{Acknowledgements:} I would like to express my gratitude to E. Koelink for discussions on the material of this paper and its presentation, and for providing some indispensible references concerning basic hypergeometric functions and their use in quantum group theory.\\

\end{document}